\documentclass[lectures]{pcms-l}
\usepackage{amsmath}
\usepackage{amssymb}
\usepackage{amsfonts}
\usepackage{amscd}
\usepackage{curves}
\usepackage{mathrsfs}
\usepackage{amsrefs}
\usepackage{latexsym}
\usepackage{graphicx} 

\newtheorem{prob}{Problem}

\newtheorem{ther}{Theorem}

\newtheorem{propo}{Proposition}

\newtheorem{lemm}{Lemma}

\newtheorem*
{AI}{The Adjunction Inequality}  
  
\theoremstyle{definition}

\newtheorem{conj}{Conjecture}[lecture]

\def \CPb {\overline{\mathbf{CP}}{}^{2}}
\def \CP {{\mathbf{CP}}^{2}}
\def \PO {{\mathbf{CP}}^{1}}
\def \R {\mathbf{R}}
\def \Z {\mathbf{Z}}
\def \Sig{\Sigma}

\def \SS {S^2\times S^2}

\def \vp {\varphi}
\def \la {\langle}
\def \ra {\rangle}
\def \a {\alpha}
\def \b {\beta}
\def \g {\gamma}

\def\L{\Lambda}
\def \lam {\lambda}
\def \Lam {\Lambda}
\def \G {\Gamma}
\def \l {\ell}
\def \o {\omega}
\def \s {\sigma}
\def \t {\tau}

\def \bd {\partial}
\def \x {\times}
\def \- {\!\smallsetminus\!}
\def \C {\subset}

\def \BC {\mathbf{C}}

\def \Q {\mathbf{Q}}
\def \bk {\bar{\k}}
\def \e {\epsilon}

\def\pt{\text{pt}}
\def \eu{{\text{e}}}

\def \ssw {\text{SW}}
\def \sw {\mathcal{SW}}

\def \DD {\Delta}

\def \Th {\Theta}

\def\so4{\text{SO(4)}}
\def\Spin4{\text{Spin(4)}}
\def\Spinc4{\text{Spin$^{\text{c}}$(4)}}
\def\spinc{spin$^{\text{c}}$ }
\def\su2{\text{SU(2)}}
\def\uo{\text{U(1)}}
\def\ut{\text{U(2)}}

\def \ch {\chi_{{}_h}}
\def \PO {{\mathbf{CP}}^{\,1}}
\def\bk{\bar{k}}

\begin{document}
\frontmatter

%% \setcounter{page}{5}

% Other information will be put here later.

%\tableofcontents

\mainmatter

% The \mainmatter section starts with the command 

% \LectureSeries[short title]{series title\author{authorname}} 

% This introduces the overall name of the particular lecture series.
% In the optional [short tile] we want a shorter name for use in the  
% headers.  Notice how the author's name is entered.

\LectureSeries[Lectures on $4$-Manifolds]
{Six Lectures on Four $4$-Manifolds \author{Ronald Fintushel and Ronald J. Stern}}

\address{Department of Mathematics, Michigan State University \newline
\hspace*{.375in}East Lansing, Michigan 48824}
\email{\rm{ronfint@math.msu.edu}}
\address{Department of Mathematics, University of California \newline
\hspace*{.375in}Irvine,  California 92697}
\email{\rm{rstern@math.uci.edu}}
\thanks{The first author was partially supported NSF Grant DMS0305818 and the second author by NSF Grant DMS0505080}


\section*{Introduction}

Despite spectacular advances in defining invariants for
simply connected smooth and symplectic 4-dimensional manifolds and the
discovery of effective surgical techniques, we
still have been unable to classify  simply connected
smooth manifolds up to diffeomorphism. In
these notes, adapted from six lectures given at the 2006 Park City Mathematics Institute Graduate Summer School on {\it Low
Dimensional Topology}, we will
review what we do and do not know about the existence and uniqueness
of smooth and symplectic structures on closed, simply connected
$4$-manifolds. We will focus on those surgical techniques that have been effective in altering smooth and symplectic structures and the Seiberg-Witten invariants that are used to distinguish them. In the last lecture we will then pose a possible classification scheme and test it on a few examples.

\lecture{How to construct $4$-manifolds } 

The classification of smooth $4$-manifolds remains a fascinating black box that has  yet to open.  You shake it, probe it, and sometimes the box mystically falls open only to find a pair of glasses and a companion black box that itself cannot be opened in any obvious way.  These newfound  glasses, lately fitted for physicists, usually provide a different vantage point from which to view the new and now even more fascinating  black box. Again, you shake it and probe it with renewed  vigor  and enthusiasm.

For example, you could look at smooth $4$-manifolds through the glasses of  a high dimensional topologist.  There it is known that every closed manifold of dimension greater than four has only finitely many smooth structures \cite{KS}.  It is also known that for $n >2$ there exist infinitely many topological manifolds that have the homotopy type of $\mathbf{CP}^{\,n}$ \cite{wall}. By contrast, in dimension four many (and perhaps all) $4$-manifolds have infinitely many distinct smooth structures \cite{KL4M}, and there exist only two topological manifolds  homotopy equivalent to $\CP$ \cite{freedman}. 

You could also look at smooth $4$-manifolds through the glasses of  a complex algebraic geometer. There it is known that except for the elliptic and Hirzebruch surfaces every simply connected complex surface has finitely many deformation types. However, as already mentioned, it is known that many (perhaps all) simply connected smooth $4$-manifolds have infinitely many distinct smooth structures. In addition, there are simply connected smooth (and irreducible) manifolds not even homeomorphic to a complex surface \cite{rat}. 

Viewing $4$-manifolds through the glasses of a symplectic geometer also misses many manifolds. While all known irreducible simply connected smooth $4$-manifolds (other than $S^4$) are homeomorphic (possibly changing orientation) to a symplectic manifold, for these examples infinitely many of the smooth structures which they support admit no underlying symplectic structure \cite{KL4M}. Thus, viewing $4$-manifolds through the glasses of either a symplectic or complex geometer will not recover the full story. However, the amalgam of ideas from complex and symplectic geometry has provided surgical techniques that are useful in the study of smooth $4$-manifolds.  

Since this is a summer school in low-dimensional topology, we will look at $4$-manifolds from the vantage point of a $2$- and $3$-manifold topologist, consider operations successful in those dimensions, and then study their analogues in dimension $4$. We then discuss their effect on the smooth invariants \textit{du jour}, the Seiberg-Witten invariants.  

In this first lecture we list the invariants that arise from algebraic topology and then those surgical constructions that have proven useful in constructing and altering smooth and symplectic structures on $4$-manifolds. We finish the lecture by describing some important manifolds that will be useful in later lectures.

\section{Algebraic Topology} 
The critical algebraic topological information for a closed, smooth, simply connected, oriented $4$-manifold $X$ is encoded in its Euler
characteristic $e(X)$, its signature $\s(X)$, and its type $t(X)$
($0$ if the intersection form of $X$ is even and $1$ if it is
odd). These invariants completely classify the homeomorphism type of
$X$ (\cites{donaldson,freedman}).  We recast these algebraic
topological invariants by defining $\ch(X)=(e(X)+\s(X))/4$ and $c(X)=3\s(X)+2e(X)$. When $X$ is a complex surface these are the holomorphic Euler characteristic and the  self-intersection of
the first Chern class of $X$, respectively. If $\overline{X}$ denotes $X$  with the opposite orientation, then $(c,\ch, t)(\overline{X})=(48\ch-5c,5\ch-\frac{c}{2},t)(X)$

\section{Techniques used for the construction of  simply connected smooth and symplectic $4$-manifolds}

We now define surgical operations that mimic constructions in lower dimensions. 

\smallskip

\noindent{\bf{Surgery on a torus.}}\; This operation is the $4$-dimensional analogue of Dehn surgery. Assume that $X$ contains a
homologically essential torus $T$ with self-intersection zero. Let
$N_T$ denote a tubular neighborhood of $T$. Deleting the interior of
$N_T$ and regluing $T^{2} \times D^{2}$ via a diffeomorphism $\vp:
\bd(T^{2} \times D^{2}) \to \partial (X- {\rm int }\ N_T)= \partial N_T$
we obtain a new manifold $X_{\vp}$, the result of surgery on $X$ along $T$. The manifold $X_{\vp}$ is determined by the homology class $\vp_*[\bd D^2]\in H_1(\bd(X\- N_T);\Z)$. Fix a basis $\{ \a, \b, [\bd D^2]\}$ for $H_1(\bd(X\- N_T);\Z)$, then there are integers $p$, $q$, $r$, such that 
$\vp_*[\bd D^2]= p\a+q\b+r[\bd D^2]$. We sometimes write $X_{\vp}=X_T(p,q,r)$. It is often the case that $X_{\vp}=X_T(p,q,r)$ only depends upon $r$, e.g. $T$ is contained in a cusp neighborhood, i.e. $\a$ and $\b$ can be chosen so that they bound vanishing cycles in $(X- {\rm int }\ N_T)$. We will sometimes refer to this process as a {\it generalized logarithmic transform} or an {\it $r$-surgery} along $T$.

If the complement of $T$ is simply connected and $t(X)=1$, then
$X_{\vp}=X_T(p,q,r)$ is homeomorphic to $X$. If the complement of $T$ is
simply connected and $t(X)=0$, then $X_{\vp}$ is homeomorphic to $X$
if $r$ is odd, otherwise $X_{\vp}$ has the same $c$ and $\ch$ but
with $t(X_{\vp})=1$. 

\smallskip

\noindent{\bf{Knot surgery.}}\; This operation is the $4$-dimensional analogue of sewing in a knot complement along a circle in a $3$-manifold. Let $X$ be a $4$-manifold which contains a
homologically essential torus $T$ of self-intersection $0$, and let
$K$ be a knot in $S^3$. Let $N(K)$ be a tubular neighborhood of $K$ in
$S^3$, and let $T\x D^2$ be a tubular neighborhood of $T$ in $X$. Then
the knot surgery manifold $X_K$ is defined by
\[ X_K = (X\- (T\x D^2))\cup (S^1\x (S^3\- N(K))\]
The two pieces are glued together in such a way that the homology
class $[{\text{pt}}\x \bd D^2]$ is identified with $[{\text{pt}}\x
\lam]$ where $\lam$ is the class of a longitude of $K$.  If the
complement of $T$ in $X$ is simply connected, then $X_{K}$ is
homeomorphic to $X$.

\smallskip

\noindent{\bf{Fiber sum.}}\; This operation is a $4$-dimensional analogue of sewing together knot complements in dimension $3$, where a knot in dimension $4$ is viewed as an embedded surface. Assume that two $4$-manifolds $X_{1}$ and
$X_{2}$ each contain an embedded genus $g$ surface $F_{j} \subset
X_{j}$ with self-intersection $0$. Identify tubular neighborhoods $N_{F_{j}}$ of $F_{j}$ with $F_{j}\times D^{2}$ and fix a diffeomorphism
$f:F_{1} \to F_{2}$. Then the fiber sum $X= X_{1}\#_{f} X_{2}$ of
$(X_{1}, F_{1})$ and $(X_{2},F_{2})$ is defined as $X_{1}\- N_{F_{1}}\cup_{\vp}X_{2}\- N_{F_{2}}$, where $\vp$ is $f \times$
(complex conjugation) on the boundary $\partial(X_{j}\-  N_{F_{j}})=F_{j}\times S^{1}$. We have
\[(c,\ch)(X_{1}\#_{f} X_{2})= (c,\ch)(X_1) +(c,\ch)(X_2)+ (8g-8,g-1)\]
Also $t(X_{1}\#_{f} X_{2})=1$ unless $F_{j}$ is characteristic in $X_{j}$, $j=0,1$.

\smallskip

\noindent{\bf{Branched covers.}}\; A smooth proper map $f:X \to Y$ is a \textit{d-fold branched covering} if away from the critical set $B \subset Y$ the restriction $f\vert X \-  f^{-1}(B):X \- f^{-1}(B) \to Y \- B$ is a covering map of degree $d$, and for $p \in f^{-1}(B)$ there is a positive integer $m$ so that the map $f$ is $(z,x) \to (z^m,x)$ in appropriate coordinate charts around $p$ and $f(p)$. The set $B$ is called the \textit{branch locus}
of the branched cover $f:X \to Y$. In the case of \textit{cyclic} branched covers, i.e. when the index-$d$ subgroup $\pi_1(X \- f^{-1}(B) ) \subset \pi_1(Y \- B)$ is determined by a surjection $\pi_1(Y \- B) \to \Z_d$, and $B$ is a smooth curve in $Y$, then $e(X)=d\, e(Y)-(d-1)\, e(B)$ and $\s(X)=d\,\s(Y)-\frac{d^2-1}{3d} B^2$, and it follows that
\[(c,\ch)(X)=d(c,\ch)(Y)-(d-1)e(B)(2,\frac14)-\frac{(d^2-1)}{3d}B^2(3,\frac14)\]

 \smallskip

\noindent{\bf{Blowup.}}\; This operation is borrowed from complex geometry. Form $X \#\CPb$, where $\CPb$ is the complex projective plane $\CP$ with the opposite orientation.

\smallskip

\noindent{\bf{Rational blowdown.}}\;  Let $C_{p}$ be the smooth $4$-manifold obtained by plumbing $p-1$ disk bundles over the $2$-sphere according to
the diagram

\begin{picture}(100,60)(-90,-25)
 \put(-12,3){\makebox(200,20)[bl]{$-(p+2)$ \hspace{6pt}
                                  $-2$ \hspace{96pt} $-2$}}
 \put(4,-25){\makebox(200,20)[tl]{$u_{0}$ \hspace{25pt}
                                  $u_{1}$ \hspace{86pt} $u_{p-2}$}}
  \multiput(10,0)(40,0){2}{\line(1,0){40}}
  \multiput(10,0)(40,0){2}{\circle*{3}}
  \multiput(100,0)(5,0){4}{\makebox(0,0){$\cdots$}}
  \put(125,0){\line(1,0){40}}
  \put(165,0){\circle*{3}}
\end{picture}

\noindent Then the classes of the $0$-sections have self-intersections $u_0^2=-(p+2)$ and $u_i^2=-2$, $i=1,\dots,p-2$. The boundary of $C_p$ is the lens space 
$L(p^2, 1-p)$ which bounds a rational ball $B_p$ with $\pi_1(B_p)=\Z_p$ and $\pi_1(\bd B_p)\to \pi_1(B_p)$ surjective. If $C_p$ is embedded in a $4$-manifold $X$ then the rational blowdown manifold $X_{(p)}$ of  \cite{rat} is obtained by replacing $C_p$ with $B_p$, i.e., $X_{(p)} = (X\- C_p) \cup B_p$. 

\smallskip

\noindent{\bf{Connected sum.}}\; Another  operation is the {\it{connected sum}} $X_1 \# X_2$ of two $4$-manifolds $X_1$ and $X_2$. We call a $4$-manifold {\it irreducible} if it cannot be represented as the connected sum of two manifolds except if one factor is a homotopy $4$-sphere. Keep in mind that we do not know if there exist smooth homotopy $4$-spheres not diffeomorphic to the standard $4$-sphere $S^{4}$ and that we have very little understanding of the uniqueness of connect sum decompositions of a reducible $4$-manifold.

\section{Some examples: Horikawa Surfaces}

Starting with the complex projective plane $\CP$ and $S^2 \times S^2$ we can now use these operations to construct many  complex surfaces. For example, blowing up $\CP$ nine times we obtain the rational elliptic surface $E(1)=\CP \#\, 9 \CPb$. The reason this is called  an elliptic surface is that there is a pencil of cubics in $\CP$, and by blowing up each of the nine base points ({\it{i.e.}} intersection points) of this family of cubics we obtain a family of tori (elliptic curves). Let $F$ be one of these smoothly embedded tori (fibers) with trivial normal bundle. We can then take the fiber sum of $n$ copies $E(1)$ along parallel copies of $F$  to obtain the elliptic surfaces $E(n)$. Note that 
\[(c,\ch,t)(E(n)) = (0, n, n \!\!\!\!\pmod{2})\]
The elliptic surface $E(n)$ is also a double cover of $\SS$
branched over
four disjoint copies of $S^2\x \{{\rm{pt}}\}$ together with
$2n$ disjoint copies of $\{{\rm{pt}}\}\x S^2$. The resultant branched
cover has $8n$ singular
points (corresponding to the double points in the branch set), whose
neighborhoods are cones
on
${{\mathbf{RP}}^{\,3}}$. These are desingularized in the usual way,
replacing their
neighborhoods with cotangent bundles of $S^2$. The result is $E(n)$. Equivalently, $E(n)$ is the $2$-fold  cover of $\SS$  branched over a complex curve of bidegree $(4,2n)$ in $\SS$. 
The horizontal and vertical fibrations of $\SS$ pull back to give (Lefschetz) fibrations of
$E(n)$ over $\PO$. A generic fiber of the vertical fibration is the
double cover of $S^2$
branched over $4$ points --- a torus. This describes an elliptic fibration on
$E(n)$. The generic fiber of the horizontal fibration is the double cover of
$S^2$ branched over $2n$ points, and this gives a genus $n-1$ fibration on
$E(n)$. This genus $n-1$ fibration has four singular fibers which are
the preimages of the
four $S^2\x \{{\rm{pt}}\}$'s in the branch set  together with the
spheres of self-intersection
$-2$ arising from desingularization.  The generic fiber
$F$ of the elliptic fibration meets a generic fiber
$\Sig_{n-1}$ of the horizontal fibration in two points,
$\Sigma_{n-1}\cdot F=2$.

Likewise, we can consider the the $2$-fold  cover of $\SS$ branched over a connected curve of bidgree $(6,2n)$ to obtain the Horikawa surfaces $H(2n-1)$. Note that \[ (c,\ch,t)(H(2n-1)) = (4n-8, 2n-1, 1)\]
Horikawa \cite{Horikawa} showed that  these complex surfaces support at most two deformation inequivalent complex structures. A complex surface in the other deformation type  is obtained as follows. Consider the Hirzebruch surfaces $F_{2n}$. These can be formed by taking the double of the euler class $2n$ disk bundle over the sphere. While they are all diffeomorphic to $\SS$, $n$ classifies the deformation type of the induced complex structure on $\SS$. Each $F_{2n}$ contains a sphere $S_+$ of self-intersection $2n$ and another sphere $S_-$  of self-intersection $-2n$ disjoint from $S_+$. Let $B_{n}$ be a curve in $F_{2n}$ with two components,  $B_{n+}$ and $S_-$, with $B_{n+}$ representing $5[S_+]$. The $2$-fold  cover of $F_{2n}$ branched over  $B_{n}$ is the complex surface $H'(4n-1)$ with \[(c,\ch)(H(4n-1)) = (8n-8, 4n-1) = (c,\ch)(H'(4n-1))\]
Horikawa has shown that $H(4n-1)$ and $H'(4n-1)$ are the only distinct deformation types of complex surfaces with $c=2\ch - 6$ \cite{Horikawa}.

Further, note that $t(H'(4n-1))= n\pmod2$, so that if $n$ is odd, $H(4n-1)$ and $H'(4n-1)$ are homeomorphic (and are not homeomorphic if $n=0\pmod2$ and hence represent two distinct deformation types). It is an open question whether for $n= 1\pmod2>1$ these surfaces are either diffeomorphic or symplectomorphic. Weak evidence that they are not diffeomorphic is provided by the $n=1$ case where $H(3)$ is, by construction, the elliptic surface $E(3)$ and is minimal. But $H'(3)$ contains a sphere of self-intersection $-1$ (the lift of $S_-$). Hence $H'(3)$ is not minimal and cannot be diffeomorphic to the elliptic surface $H(3)$.  Unfortunately for $n >1$ both $H(2n-1)$ and $H'(2n-1)$ have the same Seiberg-Witten invariants and we are unable to distinguish them.  We will discuss these surfaces in the last lecture.

Of course, this construction extends to the $2$-fold  covers $H(m,n)$ of $\SS$ branched over curves of bidegree $(2m, 2n)$. Then \[(c,\ch)(H(m,n))=(4(m-2)(n-2), (m-1)(n-1)+1)\]

It is a fun exercise to determine manifolds homeomorphic to some of the $H(m,n)$ as  covers of the Hirzebruch surfaces $F_r$ with branch set a disconnected curve in the spirit of the Horikawa surfaces that result in different deformation types. As for the Horikawa surfaces, It is an open problem to determine if  these homeomorphic complex surfaces are diffeomorphic. We will see in the last lecture, however, that they are all related by a single $\pm 1$  generalized logarithmic transform along homologically essential tori. 

\section{Complex Surfaces}

There are restrictions on $(c, \ch)(X)$ for a simply connected complex surface $X$ that is {\it minimal}, i.e. is not of the form $Y\#\CPb$. They either have $c=0$, i.e.  the elliptic surfaces $E(n)$ and their logarithmic transforms (which we will further discuss in Lecture 2), or they have $\ch =1$ and $c \ge 0$, i.e. $\CP$, $S^2 \times S^2$, or minimal complex manifolds homeomorphic to their blowups (which we will discuss in Lecture 5), or are surfaces of general type for which $2\ch -6 \le c <9\ch$ (which we will discuss in Lecture 6). If we remove the simple connectivity condition, then then the restrictions on $(c,\ch)$ remain except there are surfaces with $c=9\ch$ all of whom must be ball quotients and hence not simply connected (cf. \cites{persson,GS}). At the other extreme we have the Horikawa surfaces with $c=2\ch-6$. In Lecture 6 we will further discuss surfaces with $c=9\ch$ and $c=2\ch-6$.

\section{More Symplectic Manifolds} To motivate the utility of our list of operations, note that the elliptic surface $E(4)$ has a sphere with self-intersection $-4$. We can rationally blow down this sphere to obtain a (symplectic) manifold $E^{'}(4)$ with $(c, \ch, t)(E^{'}(4))=(1,4,1)$. More generally, it is easy to see that $E(n)$ contains the configuration $C_{n-2}$ whose lead sphere is a section of $E(n)$. This can be rationally blown down do obtain a simply-connected $4$-manifold $E^{'}(n)$ with $(c, \ch, t)(E^{'}(n))=(n-3,n,1)$, which according to the above restrictions on complex surfaces,  cannot support a complex structure. It turns out that these $4$-manifolds are also irreducible and symplectic \cites{rat, sym}. However, we first need some invariants that will detect when a smooth manifold is irreducible. 

\lecture{A User's Guide to Seiberg-Witten Theory} 

\section{The Set-up}
Let $S$ be a smooth complex surface and $T^*S$ its cotangent bundle. An important associated line bundle is the {\it{canonical bundle}} of $S$, $K_S=\det T^*S$, the determinant line bundle. (This is the top exterior power.) We have $c_1(K_S)=-c_1(TS)=-c_1(S)$. We often identify $K_S$ with $c_1(K_S)\in H^2(S;\Z)$ and then, using Poincar\'e duality, with an element of $H_2(S;\Z)$. For example, for the elliptic surface $E(n)$ with $F$ the homology class of an elliptic fiber, we have $K_{E(n)}=(n-2)F$.

An important property of $K_S$ is the adjunction formula: If C is an embedded holomorphic curve in $S$ and the genus of $C$ is $g$ then 
\[2g-2=C^2+K_S\cdot C\]

Gauge theory gives us a way to mimic this important class in the setting of oriented smooth $4$-manifolds. 
The critical analytical information for a smooth, closed,
oriented $4$-manifold $X$ is encoded in its Seiberg-Witten
invariant \cite{witten}. The goal of this lecture is to provide a `user's guide' to these invariants. 
For more detailed explanations one should see \cites{witten,KMgenus,M,N}.
Consider a smooth compact oriented $4$-manifold $X$ with tangent bundle $TX$. The choice of a Riemannian metric on X reduces the structure group of $TX$ to $\so4$, which may be equivalently taken as the structure group of $PX$, the bundle of tangent frames of $X$. $\pi_1(\so4)=\Z_2$ and its double covering group is $\Spin4\cong\su2\x\su2$. A {\it{spin structure}} on $X$ is a lift of $PX$ to a principal $\Spin4$-bundle $\tilde{P}X$ over $X$ such that in the diagram

\[ \unitlength 1cm
\begin{picture}(2,1)
\put(0,1){$\tilde{P}X\to PX$}
\put(.25,.85){\vector(2,-3){.4}}
\put(1.45,.85){\vector(-2,-3){.4}}
\put(.675,0){$X$}
\end{picture}\]

\noindent the horizontal map is a double cover on each fiber of $PX$.

A spin structure gives rise to spinor bundles $S^{\pm}=\tilde{P}X\x_{\su2}\BC^2$, where the action of $\su2$ on $\tilde{P}X$ arises from one of the two factors of $\Spin4\cong\su2\x\su2$. From the point of view of algebraic topology, one can think of a spin structure on $X$ as a lift

\[ \unitlength 1cm
\begin{picture}(2,1.25)
\put(0,0){$X \to B\so4$}
\put(1.25,.85){\vector(0,-1){.5}}
\put(.675,1){$B\Spin4$}
\multiput(.35,.35)(.1,.1){5}{\circle*{.025}}
\put(.7,.7){\vector(1,1){.2}}
\end{picture}\]

\noindent The obstruction to finding such a lift is the second Stiefel-Whitney class $w_2(X)\in H_2(X;\Z_2)$. One may alternatively think in terms of the transition functions \[ \{\vp_{i,j}:U_i\cap U_j\to \so4 \}\] of $PX$. A spin structure on $X$ consists of lifts
$\tilde{\vp}_{i,j}:U_i\cap U_j\to \Spin4$. In order to give a bundle $\tilde{P}X$, these lifts must satisfy the cocycle condition $\tilde{\vp}_{i,j}\circ\tilde{\vp}_{j,k}=\tilde{\vp}_{i,k}$.
From this point of view, $\tilde{P}X$ corresponds to an element $\tilde{\xi}$ of the \v{C}ech cohomology group $H^1(X;\Spin4)$ such that in the sequence
\[  \dots\to H^1(X; \Z_2) \xrightarrow{i_*} H^1(X; \Spin4) \xrightarrow{p_*} H^1(X;\so4) \xrightarrow{\delta} H^2(X;\Z_2)\to\dots\]
$p_*\tilde{\xi}=\xi$, the element which corresponds to $PX$. Note that $\delta\xi=w_2(X)$, affirming our comment above. Also note that if $X$ admits a spin structure (i.e. a lift of $\xi$), then all such lifts are in 1-1 correspondence with $H^1(X;\Z_2)$.
To each spin structure there is associated a Dirac operator $D:\G(S^+)\to \G(S^-)$, an elliptic operator which plays an important role in topology and geometry. 

In case $w_2(X)\ne 0$, $X$ admits no spin structure, but it can still admit a \spinc structure. A \spinc structure is given by a pair of rank 2 complex vector bundles $W^{\pm}$ over $X$ with isomorphisms $\det(W^+)=\det(W^-)=L$, some complex line bundle over $X$, so that locally $W^{\pm}=S^{\pm}\otimes L^{\frac12}$. To make sense of this, consider the transition maps $\{\vp_{i,j}: U_i\cap U_j\to\so4\}$ for $PX$. We can assume that our charts have overlaps $U_i\cap U_j$ which are contractible, so that we can always get lifts $\tilde{\vp}_{i,j}: U_i\cap U_j\to\Spin4$. However, if $w_2(X)\ne0$, we can never find lifts satisfying the cocycle condition. 

Similarly, suppose that we are given a complex line bundle $L$ with transition functions $\{ g_{i,j}: U_i\cap U_j\to\uo \}$. Locally these functions have square roots $(g_{i,j})^{\frac12}$. The obstruction to finding a system of square roots which satisfy the cocycle condition, i.e. to finding a global bundle $L^{\frac12}$ over $X$ such that 
$L^{\frac12}\otimes L^{\frac12}\cong L$ is $c_1(L)\pmod2$ in $H^2(X;\Z_2)$. Now suppose that $L$ is characteristic, i.e. that $w_2(X)=c_1(L)\pmod2$. The statement that $W^{\pm}$ should locally be $S^{\pm}\otimes L^{\frac12}$ means that the tensor products $\tilde{\vp}_{i,j}\otimes (g_{i,j})^{\frac12}$ should satisfy the cocycle condition. This function has values in $(\uo\x \su2\x \su2)/\{\pm 1\}=\Spinc4$, and the corresponding obstruction is $2w_2(X)=0$; so \spinc structures exist provided we can find characteristic line bundles $L$ over $X$. A theorem of Hirzebruch and Hopf states that these exist on any oriented $4$-manifold \cite{HH}. Spin$^{\text{c}}$ structures on $X$ are classified by lifts of $w_2(X)$ to $H^2(X;\Z)$ up to the action of $H^1(X;\Z_2)$. (Spin structures correspond to $0\in H^2(X,\Z)$ up to this action.) 

The group              
$\Spinc4\cong (\uo\x \su2\x \su2)/\{\pm 1\}$ fibers over         
$SO(4)\cong (\su2\x \su2)/\{\pm1\}$ with fiber $S^1\cong \uo$. A {\it{\spinc  structure}} $s$ on $X$ is a lift of $PX$ to a principal $\Spinc4$ bundle $\hat{P}_X$ over $X$. Since 
$\ut\cong(\uo\x \su2)/\{\pm1\}$, we get representations $s^{\pm}:\Spinc4\to \ut$,
and associated rank 2 complex vector bundles 
\[ W^{\pm}=\hat{P}_X\x_{s^{\pm}}\BC^2\]
called spinor bundles, and referred to above, and $L=\det(W^{\pm})$. We sometimes write $c_1(s)$ for $c_1(L)$.

As for ordinary spin structures, one has Clifford multiplication 
\[ c:T^*X\otimes W^{\pm}\to W^{\mp}\]
written $c(v,w)=v.w$ and satisfying $v.(v.w)=-|v|^2w$. Thus $c$ induces a map
\[ c: T^\ast X \to {\text{Hom}}(W^+, W^-)\]
A connection $A$ on $L$ together with the Levi-Civita connection on the tangent bundle of $X$ forms a connection $\nabla_A:\G(W^+)\to \G(T^\ast X\otimes W^+)$ on $W^+$. 
This connection, followed by Clifford multiplication, induces the Dirac operator
\[ D_A:\Gamma(W^+)\to\Gamma(W^-)\] 
Thus $D_A$ depends both on the connection $A$ and the Riemannian metric on $X$.  
The case where $L=\det(W^+)$ is trivial corresponds to a usual spin structure on $X$, and in this case we may choose $A$ to be the trivial connection and then $D_A=D:\G(S^+)\to \G(S^-)$, the usual Dirac operator.

Fix a \spinc  structure $s$ on $X$ with determinant line bundle $L$, and let $\mathcal {A}_L$ denote the affine space of connections on the line bundle $L$. Let $F_A\in\Omega^2(X)$ denote the curvature of a connection $A$ on $L$. The Hodge star operator acts as an involution on $\Omega^2(X)$. Its $\pm1$ eigenspaces are $\Omega^2_{\pm}(X)$, the spaces of self-dual and anti-self-dual 2-forms.  We have $F_A=F_A^++F_A^-$.
The bundle of self-dual $2$-forms $\Omega^2_+(X)$ 
is also associated to $\hat{P}X$ by 
$\Omega^2_+(X)\cong\hat{P}X\x_{\su2}{\mathfrak{su}}(2)$ where $\su2$ acts on its Lie algebra  ${\mathfrak{su}}(2)\cong\BC\oplus\R$ via the adjoint action.
The map 
\[ \BC\oplus\BC\to\BC\oplus\R \hspace{.25in} (z,w)\to(z\bar{w},|z|^2-|w|^2)\]
is $\su2$-equivariant, and so it induces a map
\[ q:\G(W^+)\to\Omega^2_+(X)\]

\section{The Equations}

Given a pair $(A,\psi) \in {\mathcal{A}}_X(L)\times\G(W^+)$, i.e. $A$ a connection in
$L=\det(W^{\pm})$  and 
$\psi$ a section of $W^+$, the Seiberg-Witten equations \cite{witten} are:
\begin{eqnarray*} &D_A\psi=0\\ &F_A^+=iq(\psi)\end{eqnarray*}
The gauge group $\text{Aut}(L)=\text{Map}(X,S^1)$ acts on the space of solutions to these equations via 
\[ g\cdot(A,\psi)=(A-g^{-1}dg,g\psi)\] and its
orbit space is the {\em Seiberg-Witten moduli space} $M_X(s)$. 

Some important features of the Seiberg-Witten equations are
\begin{enumerate}
\item If $(A,\psi)$ is a solution to the Seiberg-Witten equations with $\psi\ne0$ then its
stabilizer in $\text{Aut}(L)$ is trivial. Such solutions are called {\em irreducible}. The
stabilizer of a {\em reducible} solution $(A,0)$ consists of the constant maps in
$\text{Map}(X,S^1)$. This is a copy of $S^1$.
\medskip 
\item $(A,0)$ is a reducible solution if and only if $A$ is an anti-self-dual connection on the complex line bundle $L$ (i.e. if its curvature $F_A=F_A^-$, is anti-self-dual). If $b_X^+>0$ and $c_1(L)$ is nontorsion, a generic metric on $X$ admits
no such connections.
\medskip 
\item The formal dimension of the Seiberg-Witten moduli space is calculated by the
Atiyah-Singer theorem to be
\[ \dim M_X(s)=d(s)=\frac14(c_1(L)^2-(3\,\s(X)+2\,\eu(X))\] where $\eu(X)$ is the Euler number of $X$ and $\s(X)$ is its signature. Especially interesting is the case where $\dim M_X(s)=0$, since this is precisely the
condition for $X$ to admit an almost-complex structure with first Chern class equal to
$c_1(L)$.
\medskip 
\item An anti-self-dual 2-form $\eta$ on $X$ gives us a perturbation of the Seiberg-Witten equations:
\begin{eqnarray*} 
&D_A\psi=0\\
&F_A^+ =iq(\psi) +i\eta,
\end{eqnarray*}
and for a generic perturbation $\eta$, the corresponding moduli space of solutions  
$M_X(s,\eta)$ is an orientable manifold whose dimension is $\dim M_X(s)$, provided 
$M_X(s,\eta)$ contains at least one irreducible solution. (As in (2), if $b^+(X)>0$
and $c_1(L)\ne 0$, all
solutions will be irreducible for a generic choice of metric or perturbation $\eta$.)   
For simplicity we let the notation ignore this perturbation and write
$M_X(s)$ for $M_X(s,\eta)$. An orientation is given to $M_X(s)$ by fixing a `homology orientation' for $X$, that is, an orientation of $H^1(X)\oplus H^2_+(X)$.
\medskip 
\item There is a Lichnerowicz-type theorem, proved, as usual, with an application of the
Weitzenb\"{o}ck formula \cite{witten,KMgenus}: If $X$ carries a metric of positive scalar
curvature, then the only solutions of the Seiberg-Witten equations are reducible (of the
form $(A,0)$). Hence, if $b_X^+>0$, for a generic metric of positive scalar curvature,
$M_X(s)=\emptyset$.
\medskip 
\item For each $s$, the Seiberg-Witten moduli space $M_X(s)$ is compact.
\medskip 
\item There are only finitely many characteristic line bundles $L$ on $X$ for which both
$M_X(s)\ne\emptyset$ and $\dim M_X(s)\ge0$.
\end{enumerate}
Items (6) and (7) are also proved by using the Weitzenb\"{o}ck formula \cite{witten,KMgenus}.

In case $\dim M_X(s)=0$, items (4) and (6) imply that that generically,
$M_X(s)$ is a finite set of signed points (once a homology orientation has been chosen). In this case one defines the {\em Seiberg-Witten invariant} $\ssw_X(s)$ to be the signed count of these points. Generally, $(\mathcal{A}_L\x \G(W^+) )/\text{Aut}(L)$ is homotopy equivalent to ${\mathbf{CP}}^{\infty}\x \mathbf{T}^{b_1(X)}$, and its homology can be utilized to define $\ssw_X(s)$. The Seiberg-Witten invariant is a diffeomorphism invariant provided $b^+>1$.  (The $b^+=1$ case will be discussed in a later lecture.) 

\section{Adjunction Inequality}

One of the most important consequences of the Seiberg-Witten equations is the Adjunction Inequality.

\begin{AI}\cite{KMgenus} Suppose $b^+(X)>1$ and $\ssw_X(s)\ne 0$. Let $k$ be the Poincar\'e dual of $c_1(s)$. If $\Sig$ is an embedded closed surface in $X$ with self-intersection $\ge 0$ and
genus $g\ge 1$ then $2g-2\ge\Sig\cdot\Sig +|k\cdot\Sig|$.
\end{AI}

We often view $\ssw_X$ as a map from characteristic homology classes to the integers by defining $\ssw_X(k)=\sum_{PD(c_1(s))=k}\ssw_X(s)$. (Recall that \spinc structures on $X$ are classified by integral lifts of $w_2(X)$ up to the action of $H^1(X;\Z_2)$.  So if $H_1(X;\Z)$ has no $2$-torsion, and in particular, if $X$ is simply connected, then for each characteristic $k\in H_2(X;\Z)$ there is a unique \spinc structure $s$ satisfying $PD(c_1(s))=k$.) 

A class $k$ for which $\ssw_X(k)\ne 0$ is called a {\it{basic class}} of $X$. The adjunction inequality shows that basic classes satisfy a property similar to that of the canonical class of a complex surface. Property (7) above shows that (for $b^+(X)>1$), $X$ has at most finitely many basic classes.
\section{K\"{a}hler manifolds}

A K\"{a}hler surface is a complex surface with a metric $g$ such that $g(Jx,y)=\o(x,y)$ is a symplectic form. Each simply connected complex surface admits a K\"{a}hler structure. A K\"{a}hler surface has a distinguished \spinc structure $s_K$ with $c_1(s_K)=K_X$, the canonical class of $X$. ($K_X=-c_1(TX)$.)

\begin{ther}[\cite{witten}] If $X$ is a minimal K\"{a}hler surface with $b^+(X)>1$ then for its canonical class $|\ssw_X(\pm s_K)|=1$. Furthermore, if $c_1^2(X)>0$ then $\ssw_X(s)=0$ for all other \spinc structures. 
\end{ther}
\noindent (`Minimal' means that $X$ contains no embedded holomorphic 2-spheres with self-intersection equal to $-1$.)

Another important basic fact is that $\ssw_X(-s)=(-1)^{(e+\s)/4}\,\ssw_X(s)$.
It is useful to combine all the Seiberg-Witten invariants of $X$ into a single invariant in the integral group ring $\Z H_2(X;\Z)$, where for each $\a\in H_2(X;\Z)$ we let $t_\a$ denote the corresponding element in $\Z H_2(X;\Z)$. (Note that  $t_{\a}^{-1}=t_{-\a}$ and $t_0=1$.)
We then view the Seiberg-Witten invariant of $X$ as the Laurent polynomial
\[\sw_X = \sum \ssw_X(k)\,t_{k} \]
where the sum is taken over all characteristic elements $k$ of  $H_2(X;\Z)$. 

For example, for a minimal K\"{a}hler surface with $b^+>1$ and $c_1^2>0$ we have $\sw_X=t_K\pm t_K^{-1}$. The $K3$-surface, $E(2)$, is a K\"{a}hler surface with $b^+=3$ and $c_1=0$. Hence $\ssw_{K3}(0)=1$. Adjunction inequality arguments can be used to show that there are no other nontrivial Seiberg-Witten invariants of the $K3$-surface; so $\sw_{K3}=1$.
\section{Blowup Formula}

The blowup of a manifold $X$ is $X\#\CPb$. We may identify $H_2(X\#\CPb;\Z)$ with $H_2(X;\Z)\oplus H_2(\CPb;\Z)$, and let $E$ denote the generator of $H_2(\CPb;\Z)=\Z$ which is represented by a sphere, also called $E$, the exceptional curve. This sphere is oriented so that if $S$ is a complex surface with canonical class $K_S$, then the canonical class of its blowup is $K_S+E$. The ``Blowup Formula" relates the Seiberg-Witten invariant of $X\#\CPb$ to that of $X$. It takes an especially straightforward form when $X$ has {\it{simple type}}, that is, for each basic class $k$, the dimension, $d(k)$, of the associated Seiberg-Witten moduli space is $0$. Then
\[ \sw_{X\#\CPb}=\sw_X\cdot(e+e^{-1})\]
where $e=t_E$. The assumption of simple type is true for all symplectic manifolds with $b^+>1$ by a theorem of Taubes \cite{TGW}.

\section{Gluing Formula}
Our goal will be to use Seiberg-Witten invariants to study constructions of $4$-manifolds. The techniques will involve cutting and pasting along $3$-tori. Seiberg-Witten invariants can be defined for $4$-manifolds whose boundary is a disjoint union of $3$-tori under the extra condition that there is an $\varpi\in H^2(X;\R)$ which restricts nontrivially to each boundary component \cite{T}. In this case the invariant is an element of $\Z[[H_2(X;\Z)]]$, the ring of formal power series. For example, $\sw_{T^2\x D^2}= (t_T^{-1}-t_T)^{-1} = t_T+t_T^3+t_T^5+\cdots$. There is an important gluing theorem due to Morgan, Mrowka, and Szab\'{o} \cite{MMS}, B.D. Park \cite{P}, and in its most general form to Taubes \cite{T}:

\begin{ther}[Taubes] Suppose that $\bd X_1=\bd X_2=T^3$, and that $X=X_1\cup_{T^3} X_2$ has $b^+\ge 1$. Also suppose that there is a class $\varpi \in H^2(X;\R)$ restricting nontrivially to $H^2(T^3;\R)$. Let $j_i:X_i\to X$ be the inclusions. Then \[ \sw_X=(j_1)_*\sw_{X_1}\cdot (j_2)_*\sw_{X_2}\]
\end{ther}

When $b^+_X=1$, one gets an orientation of $H_+^2(X;\R)$ from $\varpi$: Since the restriction $i^*(\varpi)\in H^2(T^3;\R)$ is nonzero, there is a nonzero class $v\in H_2(T^3;\R)$ such that $\la i^*(\varpi),v\ra>0$. Then the condition 
$\la \a, i_*(v)\ra>0$ orients $H_+^2(X;\R)$. Now (as we shall discuss below) it makes sense to speak of $\sw^{\pm}_X$, and in Taubes' theorem, one takes $\sw^-_X$.

As an example of the use of Taubes' theorem, let $T$ be a homologically nontrivial torus of self-intersection  $0$ in $X$ with a tubular neighborhood $N_T=T\x D^2$; then 
\[ \sw_X=\sw_{X\- N_T} \cdot \frac{1}{t_T^{-1}-t_T}\]
Hence $\sw_{X\- N_T}=\sw_X\cdot (t_T^{-1}-t_T)$.

\section{Seiberg-Witten Invariants of Elliptic Surfaces}

We apply this gluing theorem to calculate the Seiberg-Witten invariants of the elliptic surfaces $E(n)$.  Recall that  $E(1)=\CP\#9\CPb$ which admits a holomorphic map to $S^2$ whose generic fiber is a self-intersection $0$ torus, $F$, and that $E(n)$ is the fiber sum $E(n-1)\#_F E(1)$. This means that $E(n)=(E(n-1)\- N_F)\cup_{T^3} (E(1)\- N_F)$. In this case, each inclusion $j_i$ is the identity. We have seen that $\sw_{E(2)}=1$. Hence
\[ 1=\sw_{E(2)}=(\sw_{E(1)\- N_F})^2\] so $\sw_{E(1)\- N_F}=\pm1$; we choose a homology orientation so that $\sw_{E(1)\- N_F}=-1$. Also $\sw_{E(2)\- N_F}\cdot\sw_{N_F}=
\sw_{E(2)}$; so $\sw_{E(2)\- N_F}=t^{-1}-t$. We then get $\sw_{E(3)}=\sw_{E(2)\- N_F}\cdot \sw_{E(1)\- N_F}= t-t^{-1}$.  Inductively, we see that $\sw_{E(n)}=(t-t^{-1})^{n-2}$, where $t=t_F$. Note that the top term $t^{n-2}$ corresponds to $K_{E(n)}=(n-2)F$.

\section{Nullhomologous tori}

Taubes' formula does not apply to neighborhoods of nullhomologous tori. We need a formula for the effect of surgery on the Seiberg-Witten invariants which applies also to surgery on nullhomologous tori. Let $T$ be a self-intersection $0$ torus embedded in $X$ with tubular neighborhood $N_T = T\x D^2$. Given a diffeomorphism $\vp:\bd(T\x D^2)\to \bd(X\- N_T)$ form $X_{\vp}=(X\- N_T)\cup_{\vp}(T\x D^2)$. The manifold $X_{\vp}$ is determined by the homology class $\vp_*[\bd D^2]\in H_1(\bd(X\- N_T);\Z)$. Fix a basis $\{ \a, \b, [\bd D^2]\}$ for $H_1(\bd(X\- N_T);\Z)$, then there are integers $p$, $q$, $r$, such that 
$\vp_*[\bd D^2]= p\a+q\b+r[\bd D^2]$. We write $X_{\vp}=X_T(p,q,r)$. (With this notation, note that $X_T(0,0,1)=X$.) We have the following important formula of Morgan, Mrowka, and Szab\'{o}:

\begin{ther} \cite{MMS} Given a class $k\in H_2(X)$:
\begin{multline*}\label{surgery formula} \sum_i\ssw_{X_T(p,q,r)}(k_{(p,q,r)}+i[T])= 
p\sum_i\ssw_{X_T(1,0,0)}(k_{(1,0,0)}+i[T]) +\\+q\sum_i\ssw_{X_T(0,1,0)}(k_{(0,1,0)}+i[T])
+r\sum_i\ssw_{X_T(0,0,1)}(k_{(0,0,1)}+i[T])
\end{multline*}
and there are no other nontrivial Seiberg-Witten invariants of $X_T(p,q,r)$.
\end{ther}

In this formula, $T$ denotes the torus $T_{(a,b,c)}$ which is the core $T^2\x {0}\C T^2\x D^2$ in each specific manifold $X_T(a,b,c)$ in the formula, and $k_{(a,b,c)}\in H_2(X_T(a,b,c))$ is any class which agrees with the restriction of $k$ in $H_2(X\- T\x D^2,\bd)$ in the diagram:

\[ \begin{array}{ccc}
H_2(X_T(a,b,c)) &\longrightarrow & H_2(X_T(a,b,c), T\x D^2)\\
&&\Big\downarrow \cong\\
&&H_2(X\- T\x D^2,\bd)\\
&&\Big\uparrow \cong\\
H_2(X)&\longrightarrow & H_2(X,T\x D^2)
\end{array}\]

\noindent Furthermore, unless the homology class $[T]$ is 2-divisible, in each term, each $i$ must be even since the classes $k_{(a,b,c)}+i[T]$ must be characteristic in $H_2(X_T(a,b,c))$.

Often this formula simplifies. For example, suppose that $\g= \vp_*[\bd D^2]$ is indivisible in $H_1(X\- N_T)$. Then there is a dual class $A\in H_3(X\- N_T, \bd)$ such that $A\cdot \g=1$. This means that $\bd A$ generates $H_2(N_{T(p,q,r)})$; so
\[  H_3(X_T(p,q,r), N_{T(p,q,r)}) \xrightarrow{\text{onto} }H_2(N_{T(p,q,r)})\xrightarrow{0}
H_2(X_T(p,q,r)) \]
So $T(p,q,r)$ is nullhomologous in $X_T(p,q,r)$. Hence in the Morgan, Mrowka, Szab\'{o} formula, the left hand side has just one term.

A second condition which simplifies the formula uses the adjunction inequality. Suppose that there is an embedded torus $\Sig$ of self-intersection $0$, such that $\Sig\cdot T(1,q,r) =1$ and $\Sig\cdot k_{(1,q,r)} = 0$. Then the adjunction inequality implies that  at most one of the classes $k_{(1,q,r)}+ i T_{(1,q,r)}$ can have $\ssw$ nonzero.

\section{Seiberg-Witten invariants for log transforms}

These gluing formulas can often be used to quickly calculate invariants,  for example, for log transforms. Consider an elliptic surface with a smooth fiber $F$ with tubular neighborhood $N_F\cong F \x D^2$. Choose a basis $\{\a,\b, [\bd D^2]\}$ for $H_1(\bd N_F;\Z)$ so that $\{\a,\b\}$ is a basis for $H_1(F;\Z)$, and hence is in the kernel of $\pi$, the projection of the elliptic fibration.  A {\it{log transform of multiplicity $r$}} is a surgery on $F$ of type $(p,q,r)$. Note that $r$, the multiplicity of the log transform is the degree of the map
\[ {\text{pt}}\x \bd D^2\xrightarrow {\vp|} \bd(X\- N_F)=\bd N_F \xrightarrow{\pi} \bd D^2\]

A {\it{vanishing cycle}} of an elliptic fibration (or for a Lefschetz fibration) is the homology class on a fiber of a loop which is collapsed in a nearby singular fiber. For example, a nodal fiber is a sphere with a transverse positive double point. Such a fiber occurs when a nonseparating loop on a smooth elliptic fiber collapses to a point. The trace of this ``collapse'' is a disk with relative self-intersection $-1$, often called a vanishing disk (or `Lefschetz thimble'). In $E(n)$, the homology of a fiber is spanned by vanishing cycles. It can be proved as a result the a log transform on $E(n)$ depends only on the multiplicity $r$.

Here is a short argument using the Morgan,-Mrowka-Szab\'{o} formula that shows why the Seiberg-Witten invariants only depend upon $r$. Suppose that $\a$ and $\b$ are the vanishing cycles. In $E(n)_F(1,0,0)$, the union of a vanishing disk bounded by $\a$ and the surgery disk form a sphere $E$ of self-intersection $-1$, and $E\cdot T_{(1,0,0)}=1$ (where $T_{(1,0,0)}$ is the core torus of the surgery). Blow down $E$ to get a $4$-manifold $Y$. (In other words, notice that the existence of $E$ forces $E(n)_F(1,0,0)$ to equal $Y\# \CPb$.) It can be easily seen that $Y$ contains a torus $\L$ of square $+1$, and in fact $T_{(1,0,0)}=\L-E$. The adjunction inequality then implies that $\sw_Y=0$; so the blowup formula says that $\sw_{E(n)_F(1,0,0)}=0$, and we see that only the ``$r"$-term of the Morgan-Mrowka-Szab\'{o} formula is nontrivial. 

We now use the Taubes formula for the full calculation. Suppose we perform a log transform of order $r$ on $E(n)$; then the result is $E(n;r)=(E(n)\- N_F)\cup_j (T^2\x D^2)$. Let $t$ be the class in the integral group ring corresponding to $T_r=T_{(p,q,r)}\in H_2(E(n;r))$. $T_r$ is a multiple torus in the sense that in $H_2(E(n;r))$, the generic fiber is $F=r\, T_r$; so $j_*(t_F)=t^r$. Hence
\begin{multline*}   \sw_{E(n;r)} = j_*(\sw_{E(n)\- N_F})\cdot\frac{1}{t-t^{-1}} = \frac{(t^r-t^{-r})^{n-1}}{t-t^{-1}}=\\ \sw_{E(n)}\cdot (t^{r-1}+t^{r-3}+\cdots+\t^{1-r}) 
\end{multline*}  

Each simply connected elliptic surface is the result of $0$, $1$, or $2$ log transforms on $E(n)$, $n\ge 1$, of relatively prime orders.  (In this notation, $r=1$ corresponds to no log transform.) Thus we complete the calculation of Seiberg-Witten invariants of simply connected elliptic surfaces with $b^+>1$ by noting that a similar argument shows (for $n>1$)
\[\sw_{E(n;r,s)} = \frac{(t^{rs}-t^{-rs})^{n}} {(t^r-t^{-r})(t^s-t^{-s})} \]

\bigskip
\lecture{Knot Surgery}
\section{The Knot Surgery Theorem} Let $X$ be a $4$-manifold containing an embedded torus $T$ of self-intersection $0$, and let $K$ be a knot in $S^3$. {\em Knot surgery} on $T$ is the result of replacing a tubular neighborhood $T\x D^2$ of $T$ with $S^1$ times the exterior $S^3\- N_K$ of the knot \cite{KL4M}:
\[ X_K = \left( X\- (T\x D^2)\right) \cup \left(S^1\x(S^3\- N_K)\right) \]
where $\bd D^2$ is identified with a longitude of $K$. This description doesn't necessarily determine $X_K$ up to diffeomorphism; however, under reasonable hypotheses, all manifolds obtained from the same $(X,T)$ and $K\C S^3$ will have the same Seiberg-Witten invariant. Knot surgery is a homological variant of torus surgery in the sense that a torus surgery is the process that removes a $T^2\x D^2$ and reglues it, whereas knot surgery removes a $T^2\x D^2$ and replaces it with a homology $T^2\x D^2$. 

Here is an alternative description of knot surgery: Consider a knot $K$ in $S^3$, and let $m$ denote a meridional circle to $K$. Let $M_K$ be the 3-manifold obtained by performing $0$-framed surgery on $K$. The effect of such a surgery is to span a longitude of $K$ with a disk.
The meridian $m$ can also be viewed as a circle in $M_K$. In $S^1\x M_K$ we
have the smooth torus $T_m=S^1\x m$ of self-intersection $0$. Let
$X_K$ denote the fiber sum
\[ X_K=X\#_{T=T_m}(S^1\x M_K)=\left( X\-  (T\x D^2)\right) \cup \left((S^1\x M_K)\- 
(T_m\x D^2)\right) \]
As above, the two pieces are glued together so as to
preserve the homology class $[{\text{pt}}\x \bd D^2]$. Because
$M_K$ has the homology of
$S^2\x S^1$ with the class of $m$ generating $H_1$, the complement $(S^1\x M_K) \-  (T\x D^2)$ has the homology of $T^2\x D^2$. Thus $X_K$ has the same
homology (and intersection pairing) as $X$. 

Let us make the additional assumption that $\pi_1(X)=1=\pi_1(X\- T)$. Then, since
the class of $m$ normally generates
$\pi_1(M_K)$; the fundamental group of $M_K\x S^1$ is normally generated by the image of $\pi_1(T)$, and it follows from Van Kampen's Theorem that $X_K$
is simply connected. Thus $X_K$ is homotopy equivalent to $X$. Also, in order to define Seiberg-Witten invariants, the oriented 4-manifold $X$ must also be
equipped with an orientation of $H^2_+(X;\R)$. The manifold $X_K$ inherits
an orientation as well as an orientation of $H^2_+(X_K;\R)$ from $X$.

For example, consider knot surgery on a fiber $F$ of the elliptic surface $E(2)$ (the $K3$-surface). Recall that $\sw_{E(2)}=1$. The elliptic fibration $E(2)\to S^2$ has a section $S$ which is a sphere of square $-2$. The homology class $S+F$ is represented by a torus of square $0$ which intersects a generic fiber once. Apply the adjunction inequality to this class to see that $mF$ cannot have a nonzero Seiberg-Witten invariant unless $m=0$. (Of course we already knew this, but the point is that it is the apparatus of the adjunction inequality that is forcing $\sw_{E(2)}=1$.) Now do knot surgery on $F$ with a knot $K$ of genus $g$. In $E(2)_K$ we no longer have the section $S$; the normal disk $D$ to $F$ has been removed from $S$. In its place there is a surface $S'$ of genus $g$ formed from $S\- D$ together with a Seifert surface of the knot $K$. The class $S'$ still has self-intersection $-2$, and $S'+F$ is a class represented by a genus $g+1$ surface of self-intersection $0$. The fiber $F$ still intersects $S'+F$ once. Apply the adjunction inequality to this class to test whether we can now have $\ssw_{E(2)_K}(mF)\ne 0$:
\[ 2(g+1)-2\ge (S'+F)\cdot (S'+F) +|mF\cdot (S'+F)| = |m| \]
Thus $mF$ has the possibility of having a nonzero Seiberg-Witten invariant if $|m| \le 2g$ (and is even since $E(2)_K$ is spin). Thus performing knot surgery gives us the possibility of constructing $4$-manifolds with interesting Seiberg-Witten invariants.
In fact:

\begin{ther}\cite{KL4M} Suppose that $b^+(X)> 1$ and $\pi_1(X)=1=\pi_1(X\- T)$ and that $T$ is a homologically essential torus of self-intersection $0$. Then $X_K$ is homeomorphic to $X$ and 
\[ \sw_{X_K} = \sw_X\cdot \DD_K(t^2) \]
where $t= t_T$ and $\DD_K$ is the symmetrized Alexander polynomial of $K$.
\end{ther}

In particular, $\sw_{E(2)_K}=\DD_K(t_F^2)$.
It was shown by Seifert that any symmetric Laurent polynomial
$ p(t)=a_0+\sum\limits_{i=1}^na_i(t^i+t^{-i})$ whose coefficient sum $p(1)=\pm1$
is the Alexander polynomial of some knot in $S^3$. It follows that the family of smooth $4$-manifolds homeomorphic to the $K3$-surface is at least as rich as this family of Alexander polynomials. Also note that since $\sw_{E(2)_K}(1)=\pm1$ and 
$\sw_{E(2;r)}(1)=r$, these manifolds are not diffeomorphic to a log transform, or any number of log transforms, of $K3$.

Note that if $\bar{K}$ is the mirror image knot to $K$ in $S^3$ then $S^1\x (S^3\- N_K)\cong S^1\x (S^3\- N_{\bar{K}})$ since we may view this construction as revolving the knot exterior about an axis. At $180^{\text{o}}$ in $S^1\x (S^3\- N_K)$ we see $S^3\- N_{\bar{K}}$. Thus $X_{\bar{K}}\cong X_K$. There are currently no other known examples of inequivalent knots which give diffeomorphic manifolds via knot surgery.
\section{Proof of Knot Surgery Theorem: the role of nullhomologous tori}
The rest of this section will be devoted to a presentation of the proof of the knot surgery theorem as given in \cite{KL4M}. This proof depends on the description of the Alexander polynomial of a knot in terms of the `knot theory macare\~na':
\[ \DD_{K_+}(t)=\DD_{K_-}(t)+(t^{1/2}-t^{-1/2})\cdot\DD_{K_0}(t)\]
where $K_+$ is an oriented knot or link, $K_-$ is the result of
changing a single oriented positive (right-handed) crossing in $K_+$ to a
negative (left-handed) crossing, and $K_0$ is the result of resolving the crossing
as shown in Figure~1.

Note that if $K_+$ is a knot, then so is $K_-$, and $K_0$ is a 2-component
link.  If $K_+$
is a 2-component link, then so is $K_-$, and $K_0$ is a knot.

\centerline{\unitlength 1cm
\begin{picture}(9,4)
\put (1,1){\vector(1,2){1.12}}
\put (2,1){\line(-1,2){0.45}}
\put (1.45,2.1){\vector(-1,2){.575}}
\put (1.25,0.5){$K_+$}
\put (4,1){\line(1,2){0.45}}
\put (4.55,2.1){\vector(1,2){.575}}
\put (5,1){\vector(-1,2){1.12}}
\put (4.25,0.5){$K_-$}
\put (7,1){\line(1,4){.255}}
\put (8,1){\line(-1,4){.255}}
\put (7.255,2){\vector(-1,4){.3}}
\put (7.745,2){\vector(1,4){.3}}
\put (7.25,0.5){$K_0$}
\put (3.8,-.15){Figure 1}
\end{picture}}

\smallskip

\centerline{\unitlength .5cm
\begin{picture}(5.5,6)
\put (2.25,4){\line(-3,-4){1.5}}
\put (2.75,4){\line(3,-4){1.5}}
\put (2.25,4.3){$K_+$}
\put (0.25,1.2){$K_-$}
\put (4,1.2){$K_0$}
\put (1.3,0){Figure 2}
\end{picture}}

It is proved in \cite{KL4M} that one can start with a knot $K$ and perform macare\~na moves so as
to build a
tree starting from $K$ and at each stage adding the bifurcation of Figure 2,
where each $K_+$, $K_-$, $K_0$ is a knot or 2-component link, and
so that at the
bottom of the tree we obtain only unknots and split links. Then, because
for an unknot
$U$ we have $\DD_U(t)=1$, and for a split link $S$ (of more than one
component) we have
$\DD_S(t)=0$, we can work backwards using the macare\~na relation to calculate
$\DD_K(t)$.

For example, we compute the Alexander polynomial of the trefoil knot:

\centerline{\unitlength .75cm
\begin{picture}(11,6)
\put (1,3){\oval(3,4)[r]}
\put (0,5){\line(3,-4){.4}}
\put (1,3.66){\line(-3,-4){1}}
\put (1,3.66){\line(-3,4){.4}}
\put (0,3.66){\line(3,4){1}}
\put (0,3.66){\line(3,-4){.4}}
\put (1,2.33){\line(-3,4){.4}}
\put (0,2.33){\line(3,-4){.4}}
\put (1,2.33){\line(-3,-4){1}}
\put (1,1){\line(-3,4){.4}}
\curve(0,5,-1,5,-1.4,4.75,-1.5,4)
\curve(-1.53,2.25,-1.55,2.25)
\curve(0,1,-1,1,-1.4,1.25,-1.5,1.75)
\curve(-1.5,1.6,-1.5,4)
\curve(-1.45,2.25,-1.43,2.24)
\put(2.5,3){\vector(0,1){.2}}
\put(.25,.5){$K_+$}
%%%%%%%%%%%%%%%%%%%%%%
\put (6,3){\oval(3,4)[r]}
\put (5,5){\line(3,-4){1}}
\put (6,3.66){\line(-3,-4){1}}
\put (6,5){\line(-3,-4){.4}}%
\put (5,3.66){\line(3,4){.4}}%
\put (5,3.66){\line(3,-4){.4}}
\put (6,2.33){\line(-3,4){.4}}
\put (5,2.33){\line(3,-4){.4}}
\put (6,2.33){\line(-3,-4){1}}
\put (6,1){\line(-3,4){.4}}
\curve(5,5,4,5,3.6,4.75,3.5,4)
\curve(3.47,2.25,3.45,2.25)
\curve(5,1,4,1,3.6,1.25,3.5,1.75)
\curve(3.5,1.6,3.5,4)
\curve(3.55,2.25,3.57,2.24)
\put(7.5,3){\vector(0,1){.2}}
\put(5.25,.5){$K_-$}
%%%%%%%%%%%
\put (11,3){\oval(3,4)[r]}
%\put (10,5){\line(3,-4){1}}
\put (11,3.66){\line(-3,-4){1}}
\put (11,5){\line(0,-1){1.33}}
\put (10,5){\line(0,-1){1.33}}
\put (10,3.66){\line(3,-4){.4}}
\put (11,2.33){\line(-3,4){.4}}
\put (10,2.33){\line(3,-4){.4}}
\put (11,2.33){\line(-3,-4){1}}
\put (11,1){\line(-3,4){.4}}
\curve(10,5,9,5,8.6,4.75,8.5,4)
\curve(8.47,2.25,8.45,2.25)
\curve(10,1,9,1,8.6,1.25,8.5,1.75)
\curve(8.5,1.6,8.5,4)
\curve(8.55,2.25,8.57,2.24)
\put(8.5,3){\vector(0,1){.2}}
\put(12.5,3){\vector(0,1){.2}}
\put(10.25,.5){$K_0$}
\end{picture}}

\noindent In the figure above, $K_+=K$ is the trefoil knot, $K_-$ is the unknot, and $K_0=H$ is the Hopf link. Thus we have \ 
$ \DD_K=1+(t^{1/2}-t^{-1/2})\cdot\DD_{H}$.
We see from the figure below that $H_-$ is the unlink and $H_0$ is the unknot; hence 
$\DD_H = 0+(t^{1/2}-t^{-1/2})\cdot 1$, and $\DD_K(t)=1+(t^{1/2}-t^{-1/2})^2 = t-1+t^{-1}$.

\centerline{\unitlength .5cm
\begin{picture}(23,6.5)
\curve(3.5,4,3.75,3.8,4,3.5,4.25,3,4.28,2.7,4.25,2.4,4,1.9,3.75,1.6,3.6,1.45,3.25,1.3,2.75,1.2,2.5,1.18,2,1.15,1.5,1.15,1,1.18,.75,1.2,.25,1.3,-.1, 1.45,-.25,1.6,-.5, 1.9,-.75, 2.4,-.78, 2.7,-.75, 3, -.5, 3.5,-.25, 3.8,0,4,.25,4.1,.75,4.2,1,4.22,1.5,4.25,2,4.25,2.5,4.22,2.75,4.2,3.25,4.1)
\curve(3.45,1.6,3.4,2.1,3.37,2.7,3.4,3.3,3.45,3.8,3.5,4.2,3.6,4.5,3.7,4.7,3.9,5,4.1,5.2,4.3,5.3,4.5,5.35,
4.7,5.35,4.9,5.3,5.1,5.2,5.3,5.05,5.5,4.7,5.6,4.5,5.7,4.2,
5.75,3.8,5.8,3.3,5.81,3.15,
5.81,2.25,5.8,2.1,5.75,1.6,5.7,1.2,5.6,.9,5.5,.7,5.3,.35,5.1,.2,4.9,.1,4.7,.05,4.5,.05,4.3,.1,4.1,.2,3.9,.4,3.7,.7,3.6,.9,3.5,1.2)
\put(5.8,2.5){\vector(0,1){.2}}
\put(4.25,2.5){\vector(0,-1){.2}}
\put(2,0.25){$H_+$}
%%%%%%%%%%%%%%%%%%%%%%%%%%%%%
\curve(12.25,4.1,12.5,4,12.75,3.8,13,3.5,13.25,3,13.28,2.7,13.25,2.4,13,1.9,12.75,1.6,12.6,1.45,12.25,1.3,11.75,1.2,11.5,1.18,11,1.15,10.5,1.15,10,1.18,9.75,1.2,9.25,1.3,8.9, 1.45,8.75,1.6,8.5, 1.9,8.25, 2.4,8.22, 2.7,8.25, 3, 8.5, 3.5,8.75, 3.8,9,4,9.25,4.1,9.75,4.2,10,4.22,10.5,4.25,11,4.25,11.5,4.22,11.75,4.2,12.25,4.1)
\curve(12.45,1.6,12.4,2.1,12.37,2.7,12.4,3.3,12.45,3.8)
\curve(12.55,4.25,12.7,4.7,12.9,5,13.1,5.2,13.3,5.3,13.5,5.35,
13.7,5.35,13.9,5.3,14.1,5.2,14.3,5.05,14.5,4.7,14.6,4.5,14.7,4.2,
14.75,3.8,14.8,3.3,14.81,3.15,
14.81,2.25,14.8,2.1,14.75,1.6,14.7,1.2,14.6,.9,14.5,.7,14.3,.35,14.1,.2,13.9,.1,13.7,.05,13.5,.05,13.3,.1,13.1,.2,12.9,.4,12.7,.7,12.6,.9,12.5,1.2)
\put(14.8,2.5){\vector(0,1){.2}}
\put(13.25,2.5){\vector(0,-1){.2}}
\put(11,0.25){$H_-$}
%%%%%%%%%%%%%%%%%%%%%%%%%%%%%%%
\curve(21.5,1.2,21.6,.9,21.7,.7,21.9,.4,22.1,.2,22.3,.1,22.5,.05,22.7,.05,22.9,.1,23.1,.2,23.3,.35,23.5,.7,23.6,.9,23.7,1.2,23.75,1.6,23.8,2.1,23.81,2.25,23.81,3.15,23.8,3.3,23.75,3.8,23.7,4.2,23.6,4.5,23.5,4.7,23.3,5.05,23.1,5.2,22.9,5.3,22.7,5.35,22.5,5.35,22.3,5.3,22.1,5.2,21.9,5,21.7,4.7,21.55,4.25,
21.75,3.8,22,3.5,22.25,3,22.28,2.7,22.25,2.4,22,1.9,21.75,1.6,21.6,1.45,21.25,1.3,20.75,1.2,20.5,1.18,20,1.15,19.5,1.15,19,1.18,18.75,1.2,18.25,1.3,17.9, 1.45,17.75,1.6,17.5, 1.9,17.25, 2.4,17.22, 2.7,17.25, 3, 17.5, 3.5,17.75, 3.8,18,4,18.25,4.1,18.75,4.2,19,4.22,19.5,4.25,20,4.25,20.5,4.22,20.75,4.2,21.25,4.1,21.45,3.8,21.4,3.3,21.37,2.72,21.4,2.1,21.45,1.6)
\put(23.8,2.5){\vector(0,1){.2}}
\put(22.25,2.5){\vector(0,-1){.2}}
\put(20,0.25){$H_0$}
\end{picture}}

\bigskip

We next need to describe a method for constructing 3-manifolds which was first
studied by W. Brakes \cite{WB} and extended by J. Hoste \cite{Hoste}. Let
$L$ be a link in $S^3$ with two {\it{oriented}}\/ components $C_1$ and
$C_2$. Fix tubular neighborhoods $N_i\cong S^1 \x D^2$ of $C_i$ with
$S^1\x({\text{pt on $\bd D^2$}})$ a longitude of $C_i$, i.e.
nullhomologous in $S^3\-  C_i$. For $n\in\Z$, let
$A_n=\bigl(\begin{smallmatrix} -1&0\\n&1 \end{smallmatrix}\bigr)$.
Note that $A_n$ takes a meridian to a meridian. We get a 3-manifold
\[ s(L;n)=(S^3\- (N_1\cup N_2))/ A_n \]
called a `sewn-up link exterior' by identifying $\bd N_1$ with $\bd N_2$ via 
a diffeomorphism inducing $A_n$ in homology.
A simple calculation shows that $H_1(s(L;n);\Z)=\Z\oplus\Z_{2\l-n}$ where
$\l\,$  is
the linking number in $S^3$ of the two components $C_1$, $C_2$, of $L$.
(See \cite{WB}.)
The second summand is generated by the meridian to either component.

J. Hoste \cite[p.357]{Hoste} has given a recipe for producing Kirby calculus
diagrams for $s(L;n)$. Consider a portion of $L$ consisting of a pair of strands, oriented in opposite directions, and separated by a band $B$ as in Figure~3. 

\begin{propo}\cite{Hoste} Let $L=C_1\cup C_2$ be an
oriented link in $S^3$. Consider a portion of $L$ consisting of a pair of strands as in Figure~3. The band sum of $C_1$ and $C_2$ is a knot $K$, and the sewn-up link exterior $s(L;n)$ is obtained from the framed
surgery on the the 2-component link on the right hand side of Figure~3. 
\end{propo}

\centerline{\unitlength 1cm
\begin{picture}(8,3.75)
\put (1.2,1){\vector(0,1){2.25}}
\put (1.9,3.25){\vector(0,-1){2.25}}
\put (1.9,1.4){\vector(0,-1){.4}}
\multiput (1.2,1.8)(0,.05){10}{\line(1,0){.7}}
\put (1.4,.6){\small{$L$}}
\put (.6,3){\small{$C_1$}}
\put (2.1,3){\small{$C_2$}}
\put (.8,1.9){\small{$B$}}
\put (3.25,2.05){$\longrightarrow$}
\put (5.2,2.3){\line(1,0){.15}}
\put (5.2,1.8){\line(1,0){.15}}
\put (5.2,1.8){\line(0,-1){1}}
\put (5.85,2.1){\oval(.75,1.75)[l]}
\put (5.6,2.3){\line(1,0){.45}}
\put (5.6,1.8){\line(1,0){.45}}
\put (6.05,2.3){\line(0,1){1}}
\put (5.2,2.3){\vector(0,1){1}}
\put (4.2,2.6){\small{$n-2\l$}}
\put (6.05,1.8){\vector(0,-1){1}}
\put (6.2,2.1){\oval(.75,1.75)[r]}
\put (6.7,2.6){\small{$0$}}
\put (5.45,.6){\small{$K$}}
\put (3,.25){Figure 3}
\end{picture}}

Now we can outline a proof of the knot surgery theorem. Begin with the resolution tree for a given oriented knot $K$. Each vertex of the tree corresponds to an oriented  knot or oriented $2$-component link. Replace each knot $K'$ in the tree with the $4$-manifold $X_{K'}$, and replace each $2$-component link $L$ with the fiber sum
\[ X_L = X\#_{T=S^1\x m} (S^1\x s(L; 2\l))\]
where $m$ is a meridian to either component. 

Suppose first that $K_-$ is a knot (and therefore so is $K_+$). We see in Figure~4 that $K_+$ is the result of $+1$ surgery on the circle $C$. The circle $C$ is nullhomologous; it bounds a punctured torus. In Figure~4 there is an obvious disk which is punctured twice by $K_-$. The punctured torus bounded by $C$ consists of this disk, punctured at the points of intersection with $K_-$ together with an annulus running `halfway around $K_-$'.  

\centerline{\unitlength 1cm
\begin{picture}(6.5,3.75)
\put (1,1){\line(1,2){.48}}
\put (1.55,2.05){\vector(1,2){.55}}
\put (.9,3.25){\vector(1,-2){1.15}}
\put (1.35,.6){\small{$K_+$}}
%%%%%%%%%%%%%%%%%%%%%%%%%%
\put (4,2){\oval(1.5,1)[l]}
\put (4,1){\line(1,2){.2}}
\put (4.3,1.6){\vector(1,2){.8}}
\put (4.8,1.4){\vector(1,-2){0.2}}
\put (4.45,2.1){\line(-1,2){.575}}
\put (5,2){\oval(1.5,1)[r]}
\put (4,1.49){\line(1,0){1}}
\put (4.35,2.49){\line(1,0){.25}}
\put (4.712,1.575){\line(-1,2){.18}}
\put (5.5,1.25){\small{+1}}
\put (4.35,.6){\small{$K_-$}}
\put (3.4,2.1){\small{$C$}}
%%%%%%%%%%%%%%%%%%%%%%%%%%%
\put (2.5,0){Figure 4}
\end{picture}}

Take the product of this with $S^1$, and glue into $X$ along $S^1\x m$ to obtain, on the one hand, $X_{K_+}$, and on the other, the result of (+1)-surgery on the nullhomologous torus $T_C=S^1\x C$. 
\[ X_{K_+} = (X_{K_-})_{T_C}(0,1,1) \equiv X_{K_-}(1)\]
where the basis for $H_1(\bd (T_C\x D^2))$ consists of $S^1\x \pt$, a pushoff of C in the punctured torus, and $\bd D^2=m_C$. The Morgan-Mrowka-Szab\'{o} formula implies that 
\[ \ssw_{X_{K_+}}(\a) = \ssw_{X_{K_-}}(\a)+\sum_i\ssw_{X_{K_-}}(0)(\a+2i[T_0])  \]
where $X_{K_-}(0) = (X_{K_-})_{T_C}(0,1,0)$ is the result of $0$-surgery on $T_C$. (Note that $T_C$ is also nullhomologous in $X_{K_+}$.) As in the concluding comments of \S1 of this lecture,
only one of the terms $\ssw_{X_{K_-}}(0)(\a+2i[T_0])$ in the sum can be nonzero. To see this, we show that there is a torus $\Lam$ of self-intersection $0$ in $X_{K_-}(0)$ such that $\Lam\cdot T_0=1$ (note $T_0=S^1\x m_C$), and such that $\Lam\cdot\a=0$ for all $\a\in H_2(X\- T)$. In fact, $\Lam$ is formed from the union of the punctured torus bounded by $C$ and the core disk of the $0$-framed surgery giving $X_{K_-}(0)$. We thus have 
\[ \sw_{X_{K_+}}=\sw_{X_{K_-}}+\sw_{X_{K_-}(0)}\]

The manifold $X_{K_-}(0)=X\#_{T=S^1\x m}(S^1\x Y)$ where $Y$ is the $3$-manifold obtained from $0$-framed surgery on both components of the link $K_-\cup C$ in $S^3$ as in Figure~5.

\centerline{\unitlength 1cm
\begin{picture}(8,3.75)
\put (1,2){\oval(1.5,1)[l]}
\put (1,1){\line(1,2){.2}}
\put (1.3,1.6){\vector(1,2){.8}}
\put (1.8,1.4){\vector(1,-2){0.2}}
\put (1.45,2.1){\line(-1,2){.575}}
\put (2,2){\oval(1.5,1)[r]}
\put (1,1.49){\line(1,0){1}}
\put (1.35,2.49){\line(1,0){.25}}
\put (1.712,1.575){\line(-1,2){.18}}
\put (2.5,1.25){\small{0}}
\put (1.3,.6){\small{$K_-$}}
\put (.5,3){\small{0}}
\put (6.25,1.7){\vector(1,2){.75}}
\put (6.75,1.7){\line(-1,2){.2}}
\put (6.4,2.25){\line(-1,2){.5}}
\put (6.25,1.5){\line(-1,-2){.35}}
\put (6.25,1.7){\line(1,0){.5}}
\put (6.25,1.5){\line(1,0){.5}}
\put (6.75,1.5){\vector(1,-2){.35}}
\put (3.5,2){$\longrightarrow$}
\put (4.75,2){${\small{s}}\ \Biggl($}
\put (7.8,2){$\Biggr)$}
\put (7.2,2) {\small{$; A_{2\l}$}}
\put (6.3,.6){\small{$K_0$}}
\put (3.3,0){Figure 5}
\end{picture}}

\vspace{.25in}

\medskip

\noindent Hoste's recipe tells us that $Y$ is the sewn-up manifold $s(K_0; 2\l)$. Hence, by definition, $X_{K_-}(0)=X_{K_0}$. We thus get 
\[ \sw_{X_{K_+}}=\sw_{X_{K_-}}+\sw_{X_{K_0}}\]

The other case to consider is where $L_-$ is an oriented $2$-component link (so also $L_+$ is a $2$-component link, and $L_0$ is a knot.). We get $L_+$ from $L_-$ by a single surgery on a  loop $U$ as in Figure~6.

\centerline{\unitlength 1cm
\begin{picture}(3.25,3.75)
\put (1,2){\oval(1.5,1)[l]}
\put (1,1){\line(1,2){.2}}
\put (1.3,1.6){\vector(1,2){.8}}
\put (1.8,1.4){\vector(1,-2){0.2}}
\put (1.45,2.1){\line(-1,2){.575}}
\put (2,2){\oval(1.5,1)[r]}
\put (1,1.49){\line(1,0){1}}
\put (1.35,2.49){\line(1,0){.25}}
\put (1.712,1.575){\line(-1,2){.18}}
\put (2.5,1.25){\small{+1}}
\put (.4,2.1){\small{$U$}}
\put (1.3,.6){\small{$L_-$}}
\put (.425,3){\small{$C_1$}}
\put (2.25,3){\small{$C_2$}}
\put (.9,0){Figure 6}
\end{picture}}

\vspace{.25in}

\medskip

Let $\l_-$ denote the linking number of the two components $C_1$ and $C_2$ of 
$L_-$. In the sewn-up manifold $s(L_-;2\l_-)$, the loop $U$ becomes nullhomologous, because according to Hoste's recipe $s(L_-;2\l_-)$ is:

\centerline{\unitlength 1cm
\begin{picture}(3,3.75)
\put (1.25,1.7){\vector(1,2){.75}}
\put (1.75,1.7){\line(-1,2){.2}}
\put (1.4,2.25){\line(-1,2){.5}}
\put (1.25,1.5){\line(-1,-2){.45}}
\put (.75,.5){\line(-1,-2){.2}}
\put (1.25,1.7){\line(1,0){.35}}
\put (1.65,1.7){\line(1,0){.1}}
\put (1.25,1.5){\line(1,0){.35}}
\put (1.65,1.5){\line(1,0){.1}}
\put (1.75,1.5){\line(1,-2){.45}}
\put (2.25,.5){\vector(1,-2){.2}}
\put(1.5,1.63){\oval(.25,.5)[r]}
\curve(1.5,1.88,1.42,1.84,1.4,1.75)
\curve(1.5,1.38,1.42,1.42,1.4,1.49)
\put(1.4,1.52){\line(0,1){.17}}
\put (2,.9){\oval(1.25,.75)[r]}
\put (1,.9){\oval(1.25,.75)[l]}
\put(.96,1.27){\line(1,0){.15}}
\put(2.04,1.27){\line(-1,0){.15}}
\put(1.2,1.27){\line(1,0){.6}}
\put(.97,.525){\line(1,0){1.1}}
\put(1.8,2.5){\small{$L_0$}}
\put(.2,.2){\small{$L_0$}}
\put(1,2.5){\small{$0$}}
\put(1.8,1.5){\small{$0$}}
\put(2.8,1){\small{$U$}}
\put (.9,-.5){Figure 7}
\end{picture}}

\vspace{.25in}

\medskip

\noindent Thus
\[ X_{L_+} = X\#_{T=S^1\x m}(S^1\x s(L_+;2\l_+)) = X_{L_-}(1) \]
where $X_{L_-}(1)$ is shorthand for the surgery manifold $(X_{L_-})_{S^1\x U}(0,1,1)$. The first equality is the definition of $X_{L_+}$ and the last equality is an exercise in Kirby calculus. Similarly we let $X_{L_-}(0)$ denote 
$(X_{L_-})_{S^1\x U}(0,1,0) = X\#_{T=S^1\x m}(S^1\x s(L_-;2\l_-)_0)$ where 
$s(L_-;2\l_-)_0$ stands for $0$ (rather than $+1$) surgery on $U$ in Figure~6.

\begin{propo} $\sw_{X_{L_-}(0)} = \sw_{X_{L_0}}\cdot (t-t^{-1})^2$ \end{propo}
\begin{proof} Cut open $s(L_-;2\l_-)$ to get $S^3\- N_{L_-}$ (where $N_{L_-}$ denotes a tubular neighborhood of $L_-$). Similarly we can cut open $s(L_-;2\l_-)_0$ to get a link exterior in $S^1\x S^2$. Take a product with $S^1$ and glue into $X$ giving $X\#_{T=S^1\x m} (S^1\x (S^1\x S^2))$ with a pair of tori of self-intersection $0$ removed. If we sew up the boundary of this manifold using the map $(1)\oplus A_{2\l_-}$, we re-obtain $X_{L_-}(0)$.

Instead, fill in the boundary components with copies of $T^2\x D^2$ to get a new manifold, $Z$. We wish to do this in such a way that when we remove a neighborhood of the new link $T^2\x \{0\} \cup T^2\x \{0\}\C Z$ and sew up the boundaries using 
$(1)\oplus A_0$, we get $X_{L_-}(0)$. (We want to be able to sew up with this particular matrix because $A_0$ identifies $S^1\x C_1$ with $S^1\x C_2$.) We can accomplish this by gluing each $T^2\x D^2$ to a boundary component using the matrix $(1)\oplus
\bigl(\begin{smallmatrix} 0&1\\ 1&\ -l_- \end{smallmatrix}\bigr)$. This matrix corresponds to $S^1\x$ (($-\l$)-framed surgery). Then, using the internal fiber sum formula of \S1, 
\[ \sw_{X\#_{T=S^1\x m}S^1\x s(L_-;2\l_-)_0} = \sw_Z\cdot (t-t^{-1})^2 \]

Now $Z=X\#_{T=S^1\x m} (S^1\x Y)$ where $Y$ is the $3$-manifold of Figure~8

\centerline{\unitlength 1cm
\begin{picture}(3.25,3.75)
\put (1,2){\oval(1.5,1)[l]}
\put (1,1){\line(1,2){.2}}
\put (1.3,1.6){\vector(1,2){.8}}
\put (1.8,1.4){\vector(1,-2){0.2}}
\put (1.45,2.1){\line(-1,2){.575}}
\put (2,2){\oval(1.5,1)[r]}
\put (1,1.49){\line(1,0){1}}
\put (1.35,2.49){\line(1,0){.25}}
\put (1.712,1.575){\line(-1,2){.18}}
\put (2.5,1.75){\small{$0$}}
\put (.4,2.1){\small{$U$}}
\put (.425,3){\small{$C_1$}}
\put (2.25,3){\small{$C_2$}}
\put (.35,.9){\small{$-\l_-$}}
\put (2.1,.9){\small{$-\l_-$}}
\put (.9,.25){Figure 8}
\end{picture}}
\noindent Now slide $C_1$ over $C_2$. We get $0$-surgery on $L_0$ together with a cancelling pair of handles. Thus $Y=M_{L_0}$, and the proposition is proved.
\end{proof}

It follows from the proposition and the Morgan-Mrowka-Szab\'{o} theorem that 
\[ \sw_{X_{L_+}}=\sw_{X_{L_-}} + \sw_{X_{L_0}}\cdot (t-t^{-1})^2 \]
We are now able to finish the proof of the knot surgery theorem. For a knot $K$, or an oriented $2$-component link $L$ 
and fixed $X$, we define a formal Laurent series $\Th$. For a knot $K$, define $\Th_K$ to be the quotient, $\Th_K=\sw_{X_K}/\sw_X$, and for a
2-component link define
$\Th_L=(t^{1/2}-t^{-1/2})^{-1}\cdot\sw_{X_L}/\sw_X$. It follows from from our calculations that in either case $\Th$ satisfies the relation
\[ \Th_{K_+} = \Th_{K_-} + (t-t^{-1})\cdot \Th_{K_0}. \]
Furthermore, for the unknot $U$, the manifold $X_U$ is just $X\#_T (S^2\x T^2)=X$, and so $\Th_U=1$. If $L$ is a $2$ component oriented split link, construct from $L$ the knots $K_+$ and $K_-$ as shown in Figure~9. Note that in this situation, $K_+=K_-$ and $K_0=L$.
It follows from $\Th_{K_+} = \Th_{K_-} + (t-t^{-1})\cdot \Th_{L} $ that $\Th_L=0$.

Subject to these initial values, the resolution tree and the macare\~na relation
determine $\DD_K(t)$ for
any knot $K$. It follows that $\Th_K$ is a Laurent polynomial in a single             
variable $t$, and $\Th_K(t)=\DD_K(t)$, completing the
proof of the knot surgery theorem.

\centerline{\unitlength 1cm
\begin{picture}(7,7)
\put (1,4){\line(1,2){1.12}}
\put (2,4){\line(-1,2){0.45}}
\put (1.45,5.1){\line(-1,2){.575}}
\put(2.12,6.24){\line(1,0){.75}}
\put(.875,6.24){\line(-1,0){.75}}
\put(2,4){\line(1,0){.87}}
\put(1,4){\line(-1,0){.87}}
\put(-.15,4.87){\framebox(.5,.5){$C_1$}}
\put(.125,6.24){\line(0,-1){.85}}
\put(.125,4){\line(0,1){.85}}
\put(2.6,4.87){\framebox(.5,.5){$C_2$}}
\put(2.87,6.24){\line(0,-1){.85}}
\put(2.87,4){\line(0,1){.85}}
\put (1.25,3.5){$K_+$}
%%%%%%%%%%%%%%%%%%%%%%%%%%%%%%%%
\put (5,4){\line(1,2){.45}}
\put (6.12,6.24){\line(-1,-2){.575}}
\put (4.875,6.24){\line(1,-2){1.12}}
\put(6.12,6.24){\line(1,0){.75}}
\put(4.875,6.24){\line(-1,0){.75}}
\put(6,4){\line(1,0){.87}}
\put(5,4){\line(-1,0){.87}}
\put(3.85,4.87){\framebox(.5,.5){$C_1$}}
\put(4.125,6.24){\line(0,-1){.85}}
\put(4.125,4){\line(0,1){.85}}
\put(6.6,4.87){\framebox(.5,.5){$C_2$}}
\put(6.87,6.24){\line(0,-1){.85}}
\put(6.87,4){\line(0,1){.85}}
\put(5.25,3.5){$K_-$}
%%%%%%%%%%%%%%%%%%%%%%%%%%%%%%%%
\put (3,1){\line(1,4){.255}}
\put (4,1){\line(-1,4){.255}}
\put (3.255,2){\line(-1,4){.3}}
\put (3.105,2.6){\vector(1,-4){.1}}
\put (3.745,2){\line(1,4){.3}}
\put (3.895,2.6){\vector(-1,-4){.1}}
\put(2.955,3.2){\line(-1,0){.75}}
\put(4.045,3.2){\line(1,0){.75}}
\put(2.955,1){\line(-1,0){.75}}
\put(4.045,1){\line(1,0){.75}}
\put(1.93,1.87){\framebox(.5,.5){$C_1$}}
\put(2.205,3.24){\line(0,-1){.85}}
\put(2.205,1){\line(0,1){.85}}
\put(4.52,1.87){\framebox(.5,.5){$C_2$}}
\put(4.795,3.24){\line(0,-1){.85}}
\put(4.795,1){\line(0,1){.85}}
\put (3.25,0.5){$L$}
\put (2.8,0){Figure 9}
\end{picture}} 

\bigskip

\lecture{Rational blowdowns }

\section{Configurations of Spheres and Associated Rational Balls}
Suppose that $X$ contains an embedded sphere of self-intersection $-1$, and let $N$ be its tubular neighborhood. Then $\bd N\cong S^3$ and $N\cup B^4\cong \CPb$ hence $X\cong Y\#\CPb$ for $Y= (X\- N)\cup B^4$. $Y$ is called the {\it{blowdown}} of $X$. Similarly, since the complement of a nonsingular quadric curve in $\CP$ is a tubular neighborhood $B_2$ of ${\mathbf{RP}}^{2}$, if $X$ contains an embedded sphere of self-intersection $-4$, its neighborhood, $C_2$ can be `rationally blown down' by replacing it with $B_2$. Note that $H_*(B_2;\Z)=H_*({\mathbf{RP}}^{2};\Z)$; so $H_*(B_2;\Q)=H_*(B^4;\Q)$; {\it{i.e.}} $B_2$ is a rational ball. Note that this process, like blowing down, reduces $b^-$ by one, while leaving $b^+$ unchanged.

More generally, as was mentioned in Lecture 1, neighborhoods $C_p$, $p=2,3,\dots$, of the configurations of embedded spheres

\begin{picture}(100,50)(-90,-25)
 \put(-12,3){\makebox(200,20)[bl]{$-(p+2)$ \hspace{6pt}
                                  $-2$ \hspace{96pt} $-2$}}
 \put(4,-25){\makebox(200,20)[tl]{$u_{0}$ \hspace{25pt}
                                  $u_{1}$ \hspace{86pt} $u_{p-2}$}}
  \multiput(10,0)(40,0){2}{\line(1,0){40}}
  \multiput(10,0)(40,0){2}{\circle*{3}}
  \multiput(100,0)(5,0){4}{\makebox(0,0){$\cdots$}}
  \put(125,0){\line(1,0){40}}
  \put(165,0){\circle*{3}}
\end{picture}

\noindent also bound rational balls and can be blown down. The boundary of $C_p$ is the lens space $L(p^2, 1-p)$ which bounds a rational ball $B_p$, and
 \[ \pi_1(L(p^2,1-p))=\Z_{p^2}\to \pi_1(B_p)=\Z_p\]
  is onto. If $C_p\C X$, then its rational blowdown $X_{(p)}=(X\- C_p)\cup B_p$.
This operation reduces $b^-$ by $p-1$ while leaving $b^+$ unchanged. It follows that $c(X_{(p)})=c(X)$ while $\ch(X_{(p)})=\ch(X)-(p-1)$. Also, if $X$ and $X\- C_p$ are simply connected, then so is $X_{(p)}$.

To see that $L(p^2,1-p)$ bounds a rational ball $B_p$ with $\pi_1(B_p)=\Z_p$ we can argue in several ways:

\smallskip

\noindent 1. Let ${\mathbf{F}}_{p-1}$, $p\ge 2$, be the simply connected ruled surface
whose negative
section $s_-$ has square $-(p-1)$. Let $s_+$ be a positive section (with square
$(p-1)$) and $f$ a fiber.  Then the homology classes
$s_++f$ and $s_-$ are represented by embedded $2$-spheres which intersect
each other
once  and have intersection matrix
\[ \begin{pmatrix} p+1& 1\\ 1 & -(p-1) \end{pmatrix} \] It follows that the
regular
neighborhood of this pair of $2$-spheres has boundary
$L(p^2,p-1)$. Its complement in ${\mathbf{F}}_{p-1}$ is the rational ball $B_p$.

\smallskip

\noindent 2. Consider the configuration of $(p-1)$ $2$-spheres

\centerline{\unitlength 1cm
\begin{picture}(5,2)
\put(.9,.7){$\bullet$}
\put(1,.8){\line(1,0){1.3}}
\put(2.2,.7){$\bullet$}
\put(2.3,.8){\line(1,0){.75}}
\put(3.3,.8){.}
\put(3.5,.8){.}
\put(3.7,.8){.}
\put(4,.8){\line(1,0){.75}}
\put(4.65,.7){$\bullet$}
\put(.6,1.1){$p+2$}
\put(2.2,1.1){$2$}
\put(4.65,1.1){$2$}
\end{picture}}

\noindent in $\#(p-1){\mathbf{CP}}^{\,2}$ where the spheres (from left to
right) represent
\[ -2h_1-h_2-\cdots- h_{p-1}, \ h_1-h_2, \ h_2-h_3, \dots ,
h_{p-2}-h_{p-1}\] where $h_i$
is the hyperplane class in the $i$\,th copy of ${\mathbf{CP}}^{\,2}$. Note that two adjacent spheres in this configuration have intersection number $-1$. This is the configuration $\bar{C}_p$, {\it{i.e.}} $C_p$ with the opposite orientation. Its
boundary is $L(p^2,p-1)$, and the classes of the
configuration span $H_2(\CP;\Q)$. The complement is the
rational ball $B_p$.

\smallskip

\noindent 3. A Kirby calculus description of $B_p$ consists of a $1$-handle and an unknotted $2$-handle that links it $p$ times and has framing $p-1$.

\begin{lemm} Any diffeomorphism of $L(p^2,1-p)$ extends over $B_p$.
\end{lemm}
\begin{proof} Bonahon has showed that $\pi_0(\text{\it{Diff}}(L(p^2,1-p))=\Z_2$. The Kirby calculus description (3) of $B_p$ gives a symmetric link in $S^3$. Rotation by $180^o$ gives the nontrivial element of  $\pi_0(\text{\it{Diff}}(L(p^2,1-p))$ and extends over $B_p$.
\end{proof}

This lemma implies that rational blowdown is well-defined. Here a a few examples:

The elliptic surface $E(4)$ has nine disjoint and homologically independent sections arising from fiber summing the nine sections of $E(1)$ which come from the exceptional curves $E_1,\dots, E_9\C E(1)=\CP\# 9\,\CPb$ in $E(4)$, the fiber sum of four copies of $E(1)$. Each of these sections is a sphere of square $-4$ which may be rationally blown down. We obtain manifolds with $\ch = \ch(E(4))=4$ and $c=c(E(4)+m=m$ where we have blown down $m$ of the sections. It can be checked that for $m=1,\dots,8$, these manifolds are simply connected.

The next example, whose discovery motivated our interest in rational blowdowns, shows
that  a logarithmic transform of order $p$ can be obtained by a sequence of
$(p-1)$ blowups and one rational blowdown of a natural embedding of the configuration $C_p$.
 
The homology class $F$ of the fiber of $E(n)$ can be represented  by a nodal curve, an
immersed sphere with one positive double point.  Blow up this double point ({\it{i.e.}} take the proper  transform of $F$) so that the class
$F-2e_1$  (where $e_1$ is the homology class of the exceptional divisor)  is
represented by an embedded sphere with square $-4$ in $E(n)\#\CPb$. This is the configuration $C_2$.
In $E(n)\#\CPb$, the exceptional divisor intersects the  $-4$-sphere in two positive
points. Blow up
one of  these points, {\it{i.e.}} again take a proper transform.  One obtains the
homology
classes $u_2=F-2e_1-e_2$ and $u_1=e_1-e_2$ which form the configuration
$C_3$ in $E(n)\#2\,\CPb$. 
Continuing in this fashion, $C_p$ naturally embeds in
$E(n)\#_{p-1}{\CPb}\subset E(n)\#_{p-1}{\CPb}$.

\begin{ther}[\cite{rat}]The rational blowdown of the above configuration
\[ C_p\subset
E(n)\#(p-1)\CPb\] is $E(n)_p$, the multiplicity $p$ log transform of $E(n)$.
\end{ther}

\section{Effect on Seiberg-Witten Invariants}

We next wish to study the effect of rational blowdowns on Seiberg-Witten invariants. Before undertaking this, we need a lemma, whose proof we leave as an exercise.

\begin{lemm} In the configuration $C_p$

\begin{picture}(100,50)(-90,-25)
 \put(-12,3){\makebox(200,20)[bl]{$-(p+2)$ \hspace{6pt}
                                  $-2$ \hspace{96pt} $-2$}}
 \put(4,-25){\makebox(200,20)[tl]{$U_{0}$ \hspace{25pt}
                                  $U_{1}$ \hspace{86pt} $U_{p-2}$}}
  \multiput(10,0)(40,0){2}{\line(1,0){40}}
  \multiput(10,0)(40,0){2}{\circle*{3}}
  \multiput(100,0)(5,0){4}{\makebox(0,0){$\cdots$}}
  \put(125,0){\line(1,0){40}}
  \put(165,0){\circle*{3}}
\end{picture}

\noindent take the meridian of $U_0$ to be the generator $\g$ of the homology of $L(p^2,1-p)=\bd C_p$. Then the meridian of $U_j$ represents $(1+j(p+1))\g$.
\end{lemm}

Suppose that we know the Seiberg-Witten invariants $\sw_X$ and wish to calculate $\sw_{X_{(p)}}$, the Seiberg-Witten invariant of the rational blowdown. To facilitate this discussion, it will be convenient to now view basic classes as characteristic cohomology classes. (Of course Poincar\'e duality allows us to switch between $H_2$ and $H^2$ at will.)

Let $k\in H^2(X_{(p)};\Z)$ be characteristic. Write $X_{(p)}=X'\cup B_p$ where $X'\cap B_p=L(p^2,1-p)=L$; so $X=X'\cup C_p$. The restriction $H^2(B_p;\Z)\to H^2(L(p^2,1-p);\Z)$ is multiplication, $\Z_p\xrightarrow {\cdot p} \Z_{p^2}$, and the restriction $k_L$ of $k$ to $L$ is in the image; so $k_L=mp$, for some integer $m$. Description (3) above of $B_p$, shows that $B_p$ has even type if and only if $p$ is odd. Since $k$ is characteristic, if $p$ is even, $m$ must be odd. If $p$ is odd, we can assume that $m$ is even by taking $-p< {\text{$m$ even}} < p$.

We claim that $k_L$ extends over $C_p$ as a characteristic cohomology class. Recall from description (2) of $B_p$ that the connected sum $(p-1)\CP=B_p\cup \bar{C}_p$. Hence $(p-1)\CPb=\bar{B}_p\cup C_p$, where the classes in $C_p$ are $u_0=-2E_1-E_2-,\dots -E_{p-1}, \ u_1=E_1-E_2,\dots , u_{p-2}=E_{p-2}-E_{p-1}$. (Again, view these as cohomology classes.) Define cohomology classes $\lam'_m\in H^2((p-1)\CPb;\Z)$ by:
\[ \lam'_m =-\sum_{i=1}^{(m+p-1)/2}E_i \ +\sum_{(m+p+1)/2}^{p-1}E_i \]
(Thus, for example, for $p=5$, $\lam'_{-4}=E_1+E_2+E_3+E_4$, $\lam'_{-2}=-E_1+E_2+E_3+E_4$, $\lam'_{0}=-E_1-E_2+E_3+E_4$, $\lam'_{2}=-E_1-E_2-E_3+E_4$, and $\lam'_{4}=-E_1-E_2-E_3-E_4$.) Finally, set $\lam_m=\lam'_m|_{C_p}\in H^2(C_p;\Z)$.

Note that the intersection vector of $E_1$ with $C_p$ is 
\[ (E_1\cdot U_0, E_1\cdot U_1,\dots, E_1\cdot U_{p-2})=(2,-1,0,\dots,0)\]
Thus, using the lemma above, $E_1|_L=(2-(p+2))\g=-p\g$. Similarly the intersection vector of $E_2$ with $C_p$ is $(1,1,-1,0,\dots,0)$; so again $E_2|_L=-p\g$. Similarly, one sees that each $E_i$, $1\le i\le p-1$ satisfies $E_i|_L=-p\g$. Thus $\lam_m$ restricts to $L$ as $mp\g\in H^2(L(p^2,1-p;\Z)$. This proves our claim that each characteristic cohomology class of $X'$ extends over $C_p$.

Let $k\in H^2(X_{(p)};\Z)$ be characteristic, and let $\tilde{k}$ be the (characteristic) class in $H^2(X;\Z)$ obtained by restricting $k$ to $X'$ and extending over $C_p$ by using the appropriate $\lam_m$. Note that  $\tilde{k}^2=k^2+\lam_m^2=k^2+(1-p)$. It follows easily that the dimensions of the corresponding Seiberg-Witten moduli spaces are equal, $d_X(\tilde{k})=d_{X_{(p)}}(k)$. Call such a $\tilde{k}$ a {\it{lift}} of $k$.

The main theorem on rational blowdowns, proved in \cite{rat}, is:

\begin{ther} Let $k$ be a characteristic (co)homology class of $X_{(p)}$, and let $\tilde{k}$ be any lift. Then $\sw_{X_{(p)}}(k)=\sw_X(\tilde{k})$.
\end{ther}

\section{Taut Configurations}

As a first example, we consider taut configurations. The configuration $C_p\C X$ is called {\it{taut}}, if for each basic class $\bk$ of $X$, $|\bk\cdot U_0|\le p$, and $\bk\cdot U_i=0$ for $i>0$. That is to say, the configuration is taut if its spheres satisfy the adjunction inequality with respect to all basic classes of $X$. Recall that, in general, we do not expect the adjunction inequality to hold for spheres.

For a taut configuration and basic class $\bk$ of $X$, $k|_{C_p}=(\bk\cdot U_0) {\text{PD$(U_0)$}}$ where ${\text{PD$(U_0)$}}\in H^2(C_p,\bd;\Z)$. Thus $\bk|_L=(\bk\cdot U_0)\in \Z_{p^2}$. Thus, to extend over $B_p$ we must have $\bk\cdot U_0=mp\in\Z_{p^2}$, and so $\bk\cdot U_0=\pm p$ or $0$, since $|\bk\cdot U_0|\le p$. It is an exercise (see \cite{rat}) to show that $\bk\cdot U_0=0$ is not a possibility. Hence for a taut configuration, precisely the basic classes $\bk$ which satisfy  $\bk\cdot U_0=\pm p$ descend.

Consider the result of rationally blowing down a section $S$ of $E(4)$. Recall that $S$ is a sphere of square $-4$. We obtain a simply connected manifold $Y(4)$ with $\ch=4$ and $c=1$. Since $\sw_{E(4)}=(t_F-t_F^{-1})^2$, the basic classes of $E(4)$ are $\pm 2\,F$ and $0$. Since $S\cdot F=1$, the configuration $C_2=\{ S\}$ is taut and the classes $\pm2\, F$ descend to give basic classes of $Y(4)$. It follows that $\sw_{Y(4)}=t+t^{-1}$. Comparing this with the blowup formula shows that $Y(4)$ is minimal; {\it{i.e.}} it contains no $2$-spheres of square $-1$. It follows from general algebraic surface theory that a minimal algebraic surface with $c>0$ satisfies the Noether inequality, $c\ge 2\ch-6$. This means that $Y(4)$ can admit no complex structure. (And with its opposite orientation, Y(4) cannot have a complex structure either. It is an easy exercise to use the adjunction inequality to see this.) In fact, a theorem of Symington \cite{sym} implies that $Y(4)$ does have a symplectic structure.

The other eight sections of $E(4)$ are disjoint from $S$ and so they give disjoint spheres of square $-4$ in $Y(-4)$. One of these gives yet another taut configuration which can be rationally blown down. This gives a simply connected manifold with $\ch=4$ and $c=2$, with $\sw=t_K+t_K^{-1}$. Note that $c=2\ch-6$ realizing equality in the Noether inequality. Minimal complex surfaces satisfying $c=2\ch-6$ are called Horikawa surfaces. These were discussed in Lecture 1, and it can be shown that the manifold obtained by rationally blowing down two sections of $E(4)$ is $H(4)$. More generally, each $E(n)$, $n\ge 4$, contains two disjoint taut configurations $C_{n-2}$ in wich the lead sphere $U_0$ is a section. Blowing down one of the configurations gives $Y(n)$, a symplectic $4$-manifold which admits no complex structure, and blowing down both gives the Horikawa surface $H(n)$.

We conclude our study of rational blowdowns by recalculating the Seiberg-Witten invariant of log transforms on $E(n)$. (In fact, in the days of Donaldson invariants, this is how the original calculation was accomplished.) For the sake of simplicity we will work out the case of a multiplicity $5$ log transform on $E(2)$. As above, we get $C_5$ embedded in $E(2)\# 4\,\CPb$, and the spheres of the configuration are $U_0=F-2e_1-e_2-e_3$, $U_1=e_1-e_2$, $U_2=e_2-e_3$, and $U_3=e_3-e_4$. Since $\sw_{E(2)}=1$, the blowup formula gives $\sw_{E(2)\#4\,\CPb}=(\e_1+\e_1^{-1})\dots (\e_4+\e_4^{-1})$ where $\e_i=t_{e_i}$. So the basic classes of $E(2)\#4\,\CPb$ are $\sum \pm e_i$.

When $C_5$ is rationally blown down, each $e_i$ descends to the class $F_5$ of the multiple fiber in $E(2)_5$; $F=5\,F_5$. Each class $\sum \pm e_i$ descends and has Seiberg-Witten invariant $1$. Thus 
\[ \sw_{E(2)_5}=t^4+t^2+1+t^{-2}+t^{-4} \]
$t=t_{F_5}$ (as it should).

\bigskip

\lecture{Manifolds with $b^+=1$ }

\section{Seiberg-Witten Invariants}
Let $X$ be a smooth simply connected oriented $4$-manifold with $b^+=1$. As simple examples we have $\CP$, $\CP\#k\,\CPb$, and $\SS$. As we mentioned in \S2, there are some complications in defining Seiberg-Witten invariants in this case.  Suppose that we are given a Riemannian metric $g$ for $X$. Recall that the input for the Seiberg-Witten equations consists of a pair of complex 2-plane bundles $W^{\pm}$ over $X$ and a complex line bundle $L$ satisfying $\det(W^{\pm})=L$. Given $(\psi, A)$ where $\psi$ is a section of $W^+$ and $A$ is a connection on $L$, we have the perturbed Seiberg-Witten equations:
\begin{eqnarray*} 
&D_A\psi=0\\
&F_A^+ =iq(\psi) +i\eta,
\end{eqnarray*}
where $\eta$ is a self-dual $2$-form on $X$. Suppose also that we are given a homology orientation for $X$, that is, a given orientation of $H^1(X;\R)\oplus H^2_+(X;\R)$. The Seiberg-Witten invariant depends on $g$ and $\eta$ as follows. 

Recall that the curvature $F_A$ represents the cohomology class $-2\pi i c_1(L)$, and that $H^2(X;\R)$ is the orthogonal direct sum of $H^+_g(X)$ and $H^-_g(X)$, the spaces of self-dual and antiself-dual harmonic $2$-forms for $g$. The Seiberg-Witten invariant for $(g,\eta)$ is well-defined provided there are no reducible solutions to the equations, {\it{i.e.}}  solutions with $\psi=0$, and hence $F_A^+ =i\eta$. For a fixed $L$ and $g$, this occurs only if $-2\pi c_1(L)^+=\eta\in H^+_g(X)$ or equivalently if $2\pi c_1(L)+\eta\in H^-_g(X)$. This happens if and only if $2\pi c_1(L)+\eta$ is orthogonal to $ H^+_g(X)$, a dimension $b^+$ subspace. Thus the condition that there should exist a reducible solution cuts out a codimension $b^+$ affine subspace of $H^2(X;\R)$.

It follows that if $b^+>0$, $\ssw_{X,g,\eta}(L)$ is well-defined for a generic $(g,\eta)$. If $b^+>1$ a generic path of $(g,\eta)$ contains no reducible solutions; and this is why the Seiberg-Witten invariant is an oriented diffeomorphism invariant in this case. (Its sign is determined by the homology orientation.)  In case $b^+=1$, there can be reducible solutions in one-parameter families. $H^-_g(X)$ has dimension $1$ and it is spanned by a single $g$ harmonic $2$-form $\o_g$ of norm $1$ agreeing with the homology orientation. This form is often called the {\it{period point}} for $(X,g)$.

The condition for admitting a reducible solution now becomes  $(2\pi c_1(L)+\eta)\cdot \o_g=0$. This defines a codimension $1$ affine subspace or {\it{wall}}. The Seiberg-Witten invariant is well-defined for families of $(g,\eta)$ satisfying the condition that any two can be connected by a path not crossing the wall. Let $H\in H_2(X;\R)$ be the Poincar\'e dual of $\o_g$.
Then there are two well-defined Seiberg-Witten invariants $\sw_{X,H}^{\pm}(L)$ defined by:
$\sw_{X,H}^+(L)=\sw_{X,g,\eta}(L)$ for any $(g,\eta)$ such that 
$(2\pi c_1(L)+\eta)\cdot\o_g>0$, and 
$\sw_{X,H}^-(L)=\sw_{X,g,\eta}(L)$ for any $(g,\eta)$ such that 
$(2\pi c_1(L)+\eta)\cdot\o_g<0$.

These invariants are related via a ``wall-crossing formula" \cite{KMgenus,LiLiu} which states that 
\[ \ssw_{X,H}^+(L)-\ssw_{X,H}^-(L)=(-1)^{1+\frac12 d(k)}\]
Thus for any characteristic homology class $k$, at least one of $\sw_{X,H}^{\pm}(k)$ is nonzero. For example, on the ray given by multiples of the fiber class $F$ of $E(1)$ (and with $t_F=t$) we have 
\[ \sw_{E(1),H}^-=\sum_{m=0}^{\infty}t^{2m+1},\ \ \sw_{E(1),H}^+=-\sum_{m=0}^{\infty}t^{-(2m+1)} \]
where $H$ is the class of a line in $\CP$ viewing $E(1)=\CP\#9\, \CPb$.
(The way to do this calculation is to note that in the ``chamber'' of the positive scalar curvature metric, the invariant is $0$, then use the wall-crossing formula.)

The small-perturbation Seiberg-Witten invariant is defined by 
\[ \ssw_{X,H}(k)=\begin{cases} \ssw_{X,H}^+(k) \ \ &{\text{if $k\cdot H>0$}}\\
\ssw_{X,H}^-(k) \ \ &{\text{if $k\cdot H<0$}} \end{cases} \]

There is again a wall-crossing formula which can be deduced from the one above.
If $H', H''$ are elements of positive square in $H_2(X;\R)$ with $H'\cdot H>0$ and $H''\cdot H>0$, then if $k\cdot H' <0$ and $k\cdot H''>0$,
\[ \sw_{X,H''}(k) - \sw_{X,H'}(k) = (-1)^{1+\frac12 d(k)}\]
Furthermore, in case $b^-\le 9$, one can see that $\ssw_{X,H}(k)$ is independent of $H$, and we denote it as $\ssw_X(k)$.

The small-perturbation invariant satisfies all the usual properties of the Seiberg-Witten invariant. In particular, there are only finitely many basic classes; so we get a Laurent polynomial $\sw_{X,H}$ or $\sw_X$, in case $b^-\le 9$. Furthermore, the invariant vanishes in the presence of a positive scalar curvature metric. So, for example, $\sw_{\CP\#m\, \CPb, H}=0$, for $H$ the class of a line in $\CP$ and if $m\le 9$,  $\sw_{\CP\#m\, \CPb}=0$.

$E(1)=\CP\#9\, \CPb$ is an elliptic surface. Consider the result of a log transform of multiplicity $p$ on $E(1)$. The log transform formula gives us (along the ray of a fiber):
\[ \sw_{E(1)_p,H}^-=\sw_{E(1),H}^- \cdot (\t^{p-1}+\t^{p-3}+\dots +\t^{1-p}) =\sum_{m=0}^{\infty}\t^{2m+1}
\]
where $\t$ corresponds to the multiple fiber, so $t = \t^p$. This means that we can't distinguish $\sw_{E(1)_p,H}^-$ from $\sw_{E(1),H}^-$ (at least along the ray under consideration). In fact, it can be shown that $E(1)_p$ is actually diffeomorphic to $E(1)$. 

The situation is quite different when we perform two log transforms. (In order that the result be simply connected, and hence homeomorphic to $E(1)$, there can be at most two log transforms, and their multiplicities $p,q$ must be relatively prime.) The result $E(1)_{p,q}$ is called a {\it{Dolgachev surface}}. A calculation of Seiberg-Witten invariants shows that its diffeomorphism type depends on the pair of relatively prime integers ${p,q}$. In fact, $E(1)_{2,3}$ and $E(1)$ were the first known pair of $4$-manifolds which are homeomorphic but not diffeomorphic. This was discovered in a ground-breaking paper of Simon Donaldson \cite{Dv} (1985). 

Knot surgery also gives examples of this phenomenon. For example, for the twist knots $K_n$:

\centerline{\unitlength 1.25cm
\begin{picture}(6.75,4)
\put (2.75,2.5){\oval(1.5,1)[l]}
\put (2.75,2.5){\oval(2.5,2)[l]}
\put (3.75,2.5){\oval(1.5,1)[r]}
\put (3.75,2.5){\oval(2.5,2)[r]}
\put (2.75,3.49){\line(1,0){.25}}
\put (2.75,2.99){\line(1,0){.25}}
\put (3.5,3.49){\line(1,0){.25}}
\put (2.95,3.175){\oval(1,.65)[rt]}
\put (3.6,3.325){\oval(1,.65)[lb]}
\put (3.5,2.99){\line(1,0){.25}}
\put (2.75,1.3){\framebox(1,1)}
\put (2.85,2){\Small{$2n-1$}}
\put (2.85,1.7){\Small{RH $\frac12$-}}
\put (2.85,1.35){\Small{twists}}
\put (1.75,.9) {\small{Figure 10: \ $K_n =$ twist knot}}
\end{picture}}
\vspace*{-.3in}\noindent we have $\DD_{K_n}(t)=nt-(2n-1)+nt^{-1}$, and one can calculate that $\sw_{E(1)_{K_n}}=-nt+nt^{-1}$.

The blowup formula shows that all these `exotic' manifolds remain exotic after connect summing with $\CPb$, however, it is still unknown whether there are examples of minimal exotic simply connected $4$-manifolds with $b^+=1$ and $b^->9$.

\section{Smooth structures on blow-ups of $\CP$}

Examples of interesting simply connected smooth $4$-manifolds with $b^+=1$ and $b^-<9$ are harder to come by.
In the late 1980's, Dieter Kotschick \cite{K} proved that the Barlow surface, which is homeomorphic to $\CP\#\,8 \CPb$, is not diffeomorphic to it. However,
in following years the subject of simply connected smooth $4$-manifolds with $b^+=1$  languished because of a lack of suitable examples.  In 2004, Jongil Park reignited interest in this topic by producing an example of a smooth simply connected $4$-manifold $P$ which is homeomorphic to  $\CP\#\,7 \CPb$ but not diffeomorphic to it. 
A year later, Stipsicz and Szab\'{o} produced a similar example with $b^-=6$ \cite{SS}.

Our goal in this lecture is to describe our technique (introduced in \cite{DN}) for producing infinite families of simply connected manifolds with $b^+=1$ and $b^-=6,7,8$. This technique was also used in \cite {PSS} (and \cite{DN}) to obtain examples with $b^-=5$.

The starting point for this construction is a specific elliptic fibration on $E(1)$. A simply connected elliptic surface is fibered over $S^2$ with smooth fiber a torus and with singular fibers. The most generic type of singular fiber is a nodal fiber (an immersed $2$-sphere with one transverse positive double point). Singular fibers are isolated; so he boundary of a small tubular neighborhood of a singular fiber is a torus bundle over the circle. It is therefore described by a diffeomorphism of $T^2$ called the {\it{monodromy}} of the singularity.
The monodromy of a nodal fiber is $D_a$, a Dehn twist around the `vanishing cycle' $a\in H_1(F;\Z)$, where $F$ is a smooth fiber of the elliptic fibration. As we discussed earlier, the vanishing cycle $a$ is represented by a nonseparating loop on the smooth fiber, and the nodal fiber is obtained by collapsing this vanishing cycle to a point to create a transverse self-intersection. The vanishing cycle  bounds a `vanishing disk', a disk of relative self-intersection $-1$ with respect to the framing of its boundary given by pushing the loop off itself on the smooth fiber. 

An $I_m$-singularity consists of a cycle of spheres $\{ A_i|i\in \Z_m\}$ of self-intersection $-2$. These spheres are arranged so that $A_i$ intersects $A_{i\pm1}$ once, positively and transversely, and $A_i\cap A_j=\emptyset$ for $j\ne i, \ i\pm1$. The monodromy of an $I_m$-fiber is $D_a^m$.
In particular, an $I_2$ fiber consists of a pair of $2$-spheres of self-intersection $-2$ which are plumbed at two points. The monodromy of an $I_2$-fiber is $D_a^2$, which is also the monodromy of a pair of nodal fibers with the same vanishing cycle. This means that an elliptic fibration which contains an $I_2$-fiber can be perturbed to contain two nodal fibers with the same vanishing cycle.  

The elliptic fibration on $E(1)$ which we use has an $I_8$-fiber, and $I_2$-fiber, and two other nodal fibers. (The existence of such a fibration can be found, for example, in 
\cite{Pn}.) A {\it{double node neighborhood}} $D$ is a fibered neighborhood of an elliptic fibration which contains exactly two nodal fibers with the same vanishing cycle. If $F$ is a smooth fiber of $D$, there is a vanishing class $a$ that bounds vanishing disks in the two different nodal fibers, and these give rise to a sphere $V$ of self-intersection $-2$ in $D$. It follows from our comments about $I_2$-singularities that we have a double node neighborhood $D$ in our copy of $E(1)$. 

Perform knot surgery inside the neighborhood $D$ using the twist knot $K_n$ of Figure 10. These twist knots have genus one, and an obvious genus one Seifert surface $\Sig$ can be seen in Figure 10. Consider the loop $\G$ on $\Sig$ which runs once through the clasp on top. This loop satisfies the following two properties:
\begin{itemize}
\item[(i)] $\G$ bounds a disk in $S^3$ which intersects $K$ in exactly two points.
\item[(ii)] The linking number in $S^3$ of $\G$ with its pushoff on $\Sig$ is $+1$.
\end{itemize}   
It follows from these properties that $\G$ bounds a punctured torus in $S^3\- K$. 

In the knot surgery construction, one has a certain freedom in choosing the gluing of $S^1\x (S^3\- N(K))$ to  $D\- N(F)$. We are free to make any choice as long as a longitude of $K$ is sent to the boundary circle of a normal disk to $F$. We choose the gluing so that the class of a meridian $m$ of $K$ is sent to the class of $a \x \{pt\}$ in $H_1(\bd(D\- N(F));\Z)=H_1(F\x \bd D^2;\Z)$, where $a$ is the vanishing cycle for the double node neighborhood. 

Before knot surgery, the fibration of $D$ has a section which is a disk. The result of knot surgery is to remove a smaller disk in this section and to replace it with the Seifert surface of $K_n$. Call the resulting relative homology class in $H_2(D_{K_n},\bd;\Z)$ the {\it{pseudosection}}. Thus in $D_{K_n}$, our knot-surgered neighborhood, there is a punctured torus $\Sig$ representing the pseudosection. The loop $\G$ sits on $\Sig$ and by (i) it bounds a twice-punctured disk $\DD$ in $\{pt\}\x \bd(S^3\- N(K))$ where $\bd\DD=\G\cup m_1\cup m_2$ where the $m_i$ are meridians of $K$. The meridians $m_i$ bound disjoint vanishing disks $\DD_i$ in $D\- N(F)$ since they are identified with disjoint loops each of which represents the class of $a \x \{pt\}$ in $H_1(\bd(D\- N(F));\Z)$. Hence in $D_K$ the loop $\G\C \Sig$ bounds a disk $U=\DD\cup \DD_1\cup \DD_2$. By construction, the self-intersection number of $U$ relative to the framing given by the pushoff of $\G$ in $\Sig$ is $+1-1-1=-1$. (This uses (ii).)

If we had an embedded disk $V$ in $D_K$ with boundary $\G$, with $V\cap \Sig=\G$, and with self-intersection number $0$ with respect to the framing given by $\Sig$, we could perform ambient surgery on $\Sig\C D_{K_n}$, surgering it to a $2$-sphere. Instead of $V$, we actually have an embedded disk $U$ such that  $U\cap \Sig=\G$, and with self-intersection number $-1$ with respect to the framing given by $\Sig$. We can perform ambient surgery to turn $\Sig$ in to a sphere, but now it is immersed with one positive double point. Viewing $D_{K_n}\C E(1)_{K_n}$, we get a pseudosection $S'$ of $E(1)_{K_n}$ and it is represented by an immersed sphere of square $-1$. We are then in the situation of Figure 11 where we have the configuration consisting of the immersed $2$-sphere $S'$ with a pair of disjoint nodal fibers, each intersecting $S'$ once transversely. Also, $S'$ intersects the $I_8$-fiber transversely in one point. 

\vspace{.2in}
\hspace{.75in}\includegraphics[scale=0.25]{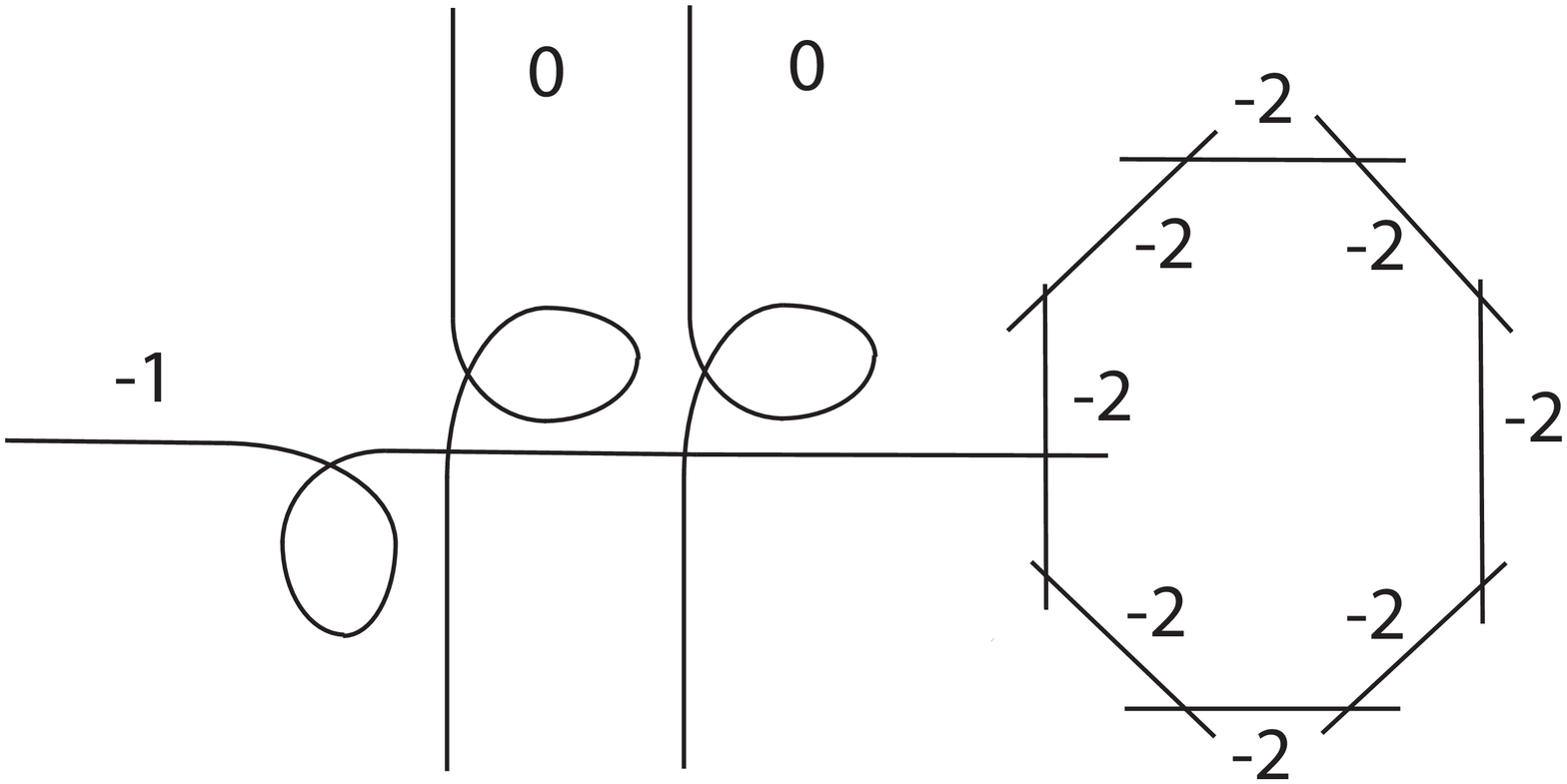}\\
\centerline{\small{Figure 11}}

\vspace{.1in} At this stage we have three possibilities:

\noindent 1.\  If we blow up the double point of $S'$, then in $E(1)_{K_n}\#\CPb$ we obtain a configuration consisting of the total transform $S''$ of $S'$, which is a sphere of self-intersection $-5$, and the sphere of self-intersection $-2$ at which $S'$ intersects $I_8$.  (See Figure 12.) This is the configuration $C_3$ which can be rationally blown down to obtain a manifold $Y_n$ with $b^+=1$ and $b^-=8$. It is easy to see that $Y_n$ is simply connected; so it is homeomorphic to $\CP\#8\,\CPb$.

\hspace{.75in}\includegraphics[scale=0.25]{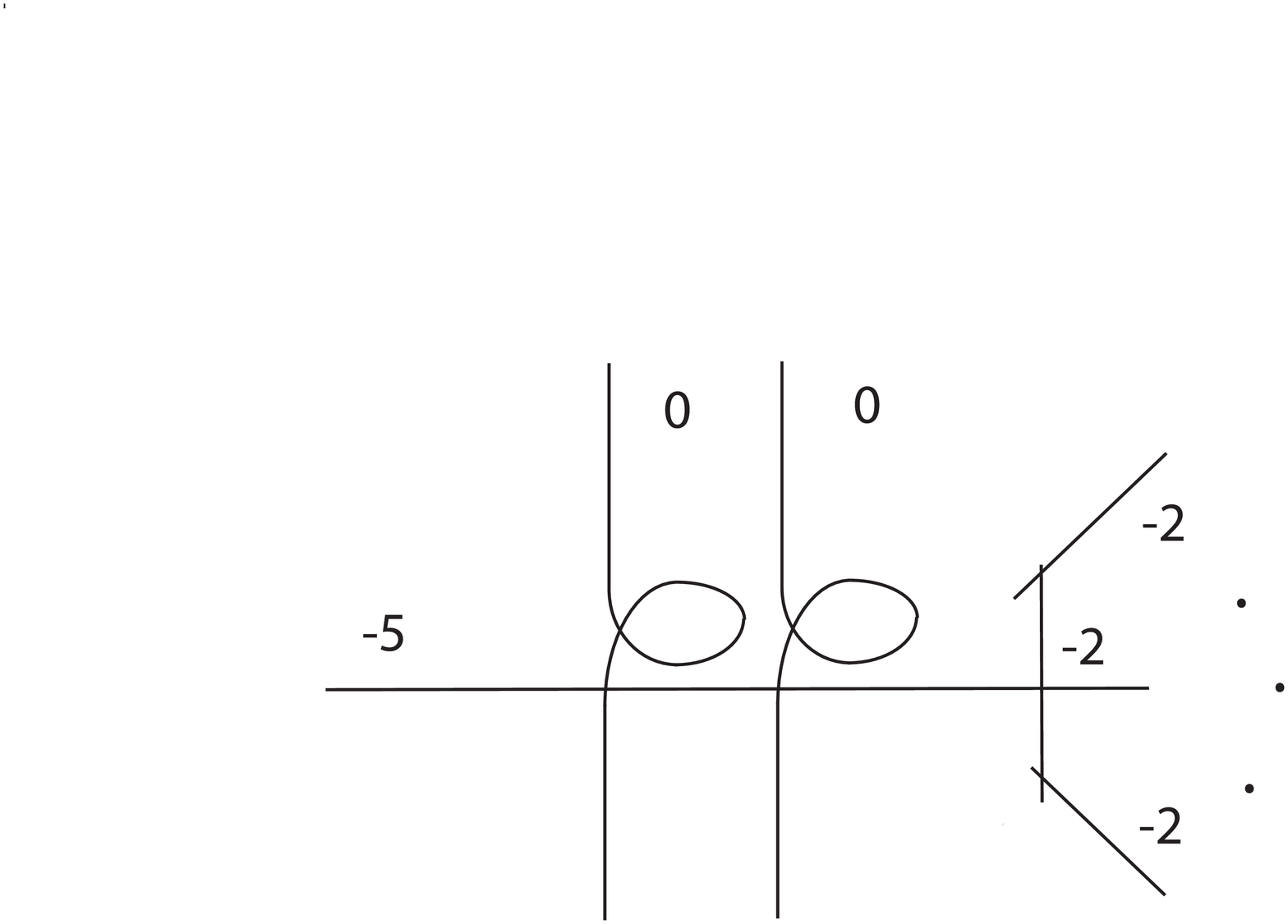}\\
\centerline{\small{Figure 12}}

\vspace{.1in}\noindent 2.\   If we blow up at the double point of $S'$ as well as at the double point of one of the nodal fibers, in $E(1)_{K_n}\#2\, \CPb$ we get a configuration of $2$-spheres consisting of $S''$, a transverse sphere $F'$ of self-intersection $-4$, and three spheres from the $I_8$-fiber.  (See Figure 13.) We can smooth the intersection of  $S''$ and $F'$ by replacing a neighborhood of the intersection point (a cone on the Hopf link) with an annulus (a Seifert surface for the Hopf link). This gives a sphere of self-intersection $-7$, and we obtain the configuration $C_5$. Rationally blowing down $C_5$ gives a manifold $X_n$ homeomorphic to $\CP\#7\,\CPb$.

\vspace{.2in}
\hspace{.75in}\includegraphics[scale=0.25]{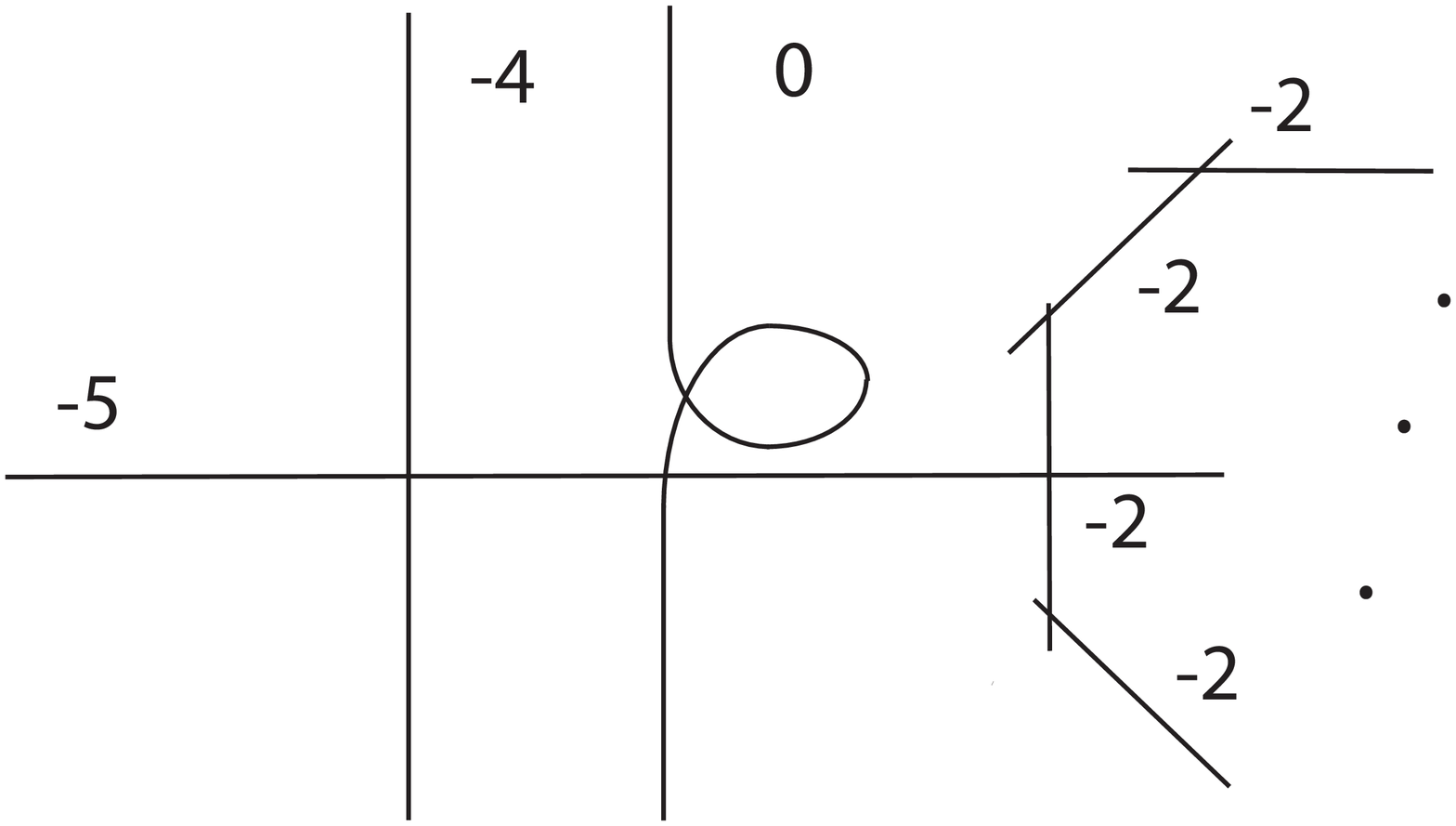}\\
\centerline{\small{Figure 13}}

\vspace{.1in}\noindent 3.\ If we blow up at the double point of $S'$ as well as at the double points of both nodal fibers, then in $E(1)_{K_n}\#3\, \CPb$ we get a configuration of $2$-spheres consisting of $S''$, two disjoint transverse spheres $F'$, $F''$ of self-intersection $-4$, and five spheres from the $I_8$-fiber. (See Figure 14.) Smoothing the intersections of  $S''$, $F'$ and $F''$ as above gives a sphere of self-intersection $-9$ and we obtain the configuration $C_7$. Rationally blowing down $C_7$ gives a manifold $Z_n$ homeomorphic to  $\CP\#6\,\CPb$.

\vspace{.2in}
\hspace{.75in}\includegraphics[scale=0.25]{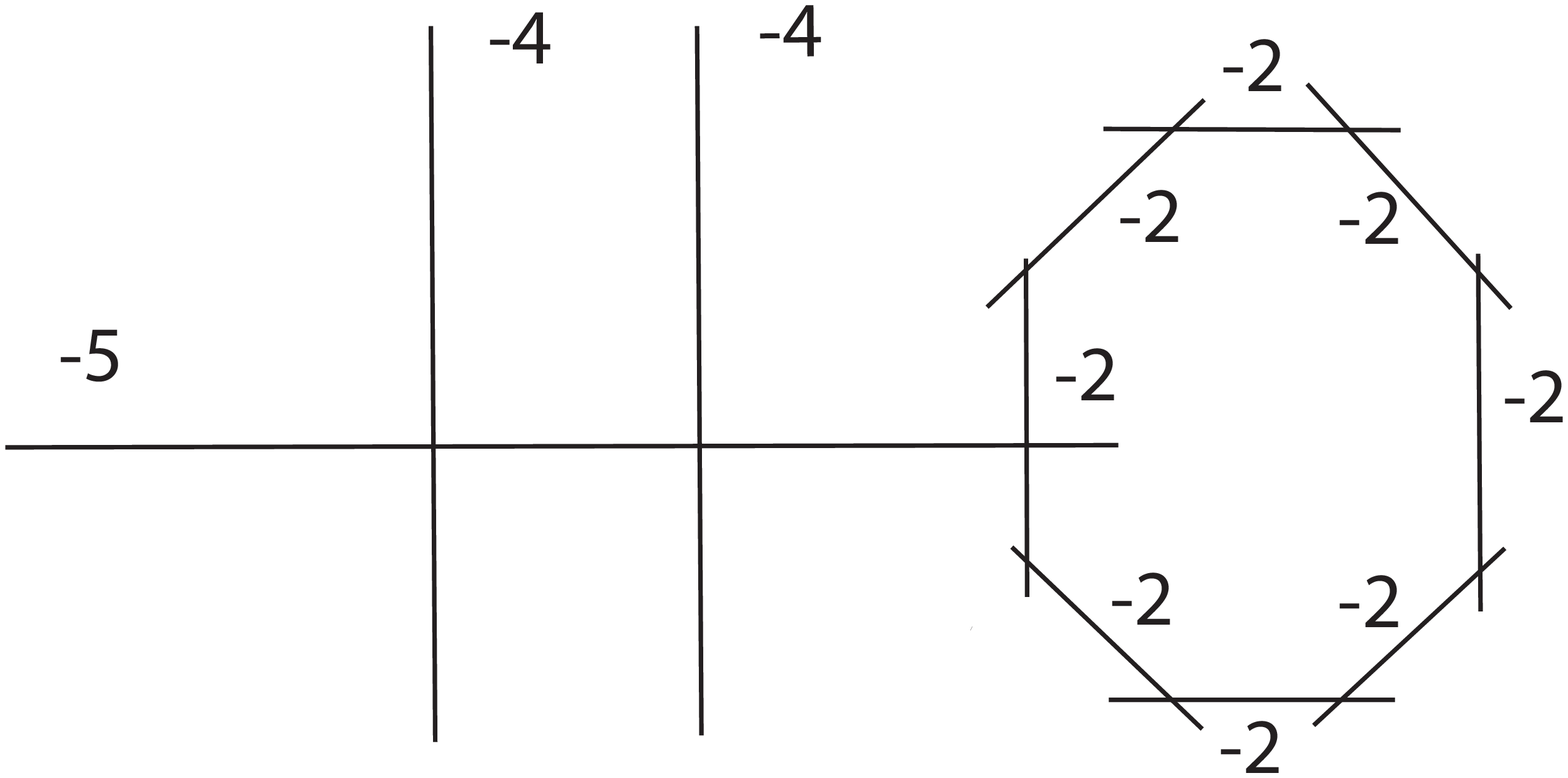}\\
\centerline{\small{Figure 14}}

\begin{ther} No two of the manifolds $Y_n$, $X_n$, $Z_n$ are diffeomorphic, and for $n\ge 2$ they are minimal. In particular, there are infinite families of mutually nondiffeomorphic $4$-manifolds homeomorphic to $CP\# k\, \CPb$, for $k=6,7,8$.
\end{ther}

We will discuss the k=7 case. The other cases are similar. Since $\sw_{E(1)_{K_n}}=-nt+nt^{-1}$, the blowup formula easily calculates $\sw_{E(1)_{K_n}\#2\,\CPb}$. If we let $\pm T$ denote the basic classes of $E(1)_{K_n}$, then $|\ssw_{E(1)_{K_n}\#2\,\CPb,h}(\pm T\pm E_1\pm E_2)| = n$,   and $\ssw_{Z_n,h}(L)=0$ for all other classes $L$. Here $h$ is the class of square $+1$ obtained as follows.

In $E(1)$ a sphere $\L$ representing $H$ intersects the fiber $F$ in 3 points. After knot surgery, this gives rise to a class $h$ of genus 3 that has $h^2=1$ and $h\cdot T =3$. (The three normal disks to a fiber that lie in the sphere $\L$ are replaced by genus one Seifert surfaces of $K_n)$.) Since the Seiberg-Witten invariant of $E(1)_{K_n}$ is well-defined, $\ssw_{E(1)_{K_n},h}(L)=\ssw_{E(1)_{K_n}}(L)$ for all characteristic $L\in H_2(E(1)_{K_n};\Z)$. The upshot of this and the wall-crossing formula is the paragraph above. 

One can now apply the criterion of \S4 to see that of the classes $\pm T\pm E_1\pm E_2$ only $\pm (T+E_1+E_2)$ descends to the rational blowdown $X_n$. So these give classes $\pm K_n\in H_2(X_n;\Z)$, and $|\ssw_{X_n}(\pm K_n)|=n$. (See \cite{DN} for details.) For all other characteristic homology classes $L$ of $X_n$, our rational blowdown theorem tells us that $\ssw_{X_n}(L) = 0,\pm1$. This proves that no two $X_m$, $X_n$, $m\ne n$, are diffeomorphic. It is an exercise in the use of the blowup formula to see that $X_n$ is minimal for $n\ge 2$.

From this vantage point it is now relatively easy to construct infinitely many distinct smooth structures on $\CP\#5\,\CPb$, as is done in \cite{DN}. We outline the construction. It begins by noting that there is an elliptic fibration on $E(1)$ containing two double node neighborhoods and one $I_6$ singularity.   Perform a knot surgery in each double node neighborhood; say using the twist knot $T(1)$ in one neighborhood and $T(n)$ in the other. In each neighborhood we need to be careful to perform the knot surgery so that the meridian of the knot is identified with the vanishing cycle. Let $V_n$ be the resultant manifold. It is an exercise to show that $V_n$ is simply-connected. 
 
Using \cite{KL4M}, one computes the Seiberg-Witten invariants of $V_n$: Up to sign the only classes in $H_2(V_n;\Z)$ with nontrivial Seiberg-Witten invariants are $T$ and $3T$, and $|\ssw_{V_n}(\pm 3T)|=n$, $|\ssw_{V_n}(\pm T)|=2n-1$.

To finish the construction, use the two double node neighborhoods to get a representative of the pseudo-section of $V_n$ which is an immersed sphere with two positive double points. If we blow up these two points we get an embedded sphere $S$ of self-intersection $-9$ in $W_n=V_n\# 2\,\CPb$. By removing one of the spheres in the $I_6$ configuration, we get $5$ spheres which, along with $S$ form the configuration $C_7$. As above, we can see that the only classes $L\in H_2(W_n;\Z)$ which have $\ssw_{W_n,h}(L)\ne 0$ and whose restriction to $C_7$ has self-intersection 
$-6$ are $L=\pm(3T+E_0+E_1)$. If we let $Q_n$ denote the result of rationally blowing down $C_7$ in $W_n$, then an argument exactly as above shows that the manifolds $Q_n$ ($n>0$) are all homeomorphic to $\CP\# 5\,\CPb$, and no two of these manifolds are diffeomorphic.
Furthermore, for $n\ge 2$, the $Q_n$ are all minimal.

\bigskip

\lecture{Putting it all together: The Geography and Botany of $4$-manifolds} 

\section {Existence: The Geography Problem} The existence part of a classification scheme for irreducible smooth (or symplectic) $4$-manifolds could take the form of determining which $(\ch, c, t) \in \Z \times \Z \times \Z_{2}$ can occur as $(\ch,c, t)(X)$ for some smooth (symplectic) $4$-manifold $X$. This is referred to as {\it the geography problem}.  Our current understanding of the geography problem is given by Figure 15 where all known simply connected smooth irreducible $4$-manifolds with nontrivial Seiberg-Witten invariants are plotted as lattice points in the $(\ch, c)$-plane. 

\setlength{\unitlength}{1in}
\begin{picture}(4.5,4.5)
\put(.2,.5){\vector(0,1){3.5}}
\put(.2,.5){\vector(1,0){4}}
\put(.2,.5){\line(1,3){1.15}}
\put(.2,.5){\line(1,2){1.7}}
\put(.8,.5){\line(3,2){3.4}}
\put(.8,.5){\line(3,1){3.4}}
\put(0,3.8){$c$}
\put(4.3,.35){$\ch$}
\put(.3,.25){Elliptic Surfaces $E(n)$  ($(\ch,c)=(n,0)$)}
\put(3.1,.3){$ c < 0$ (unknown)}
\put(2.0,2.7){$2\ch - 6 \le c \le 9\ch$}
\put(1.5,2.35){surfaces of general type}
\put(1.1,4.1){$c=9\ch$}
\put(.3,3.1){$c>9\ch$}
\put(.2,2.9){(unknown)}
\put(1.8,4.1){$c=8\ch$}
\put(1.82,3.95){$\sigma =0$}
\put(3.9,2.8){$c=2\ch - 6$}
\put(3.9,1.41){$c=\ch - 3$}
\put(1.15,3.1){$\sigma >0$}
\put(2.25,3.1){$\sigma <0$}
\put(3.10,1.9){$\ch - 3 \le c \le 2\ch -6$}
\put(2.65,1.65){symplectic with one basic class}
\put(3.0,1.50){(cf. \cite{FPS})}
\put(3.0,1.05){$0 \le c \le \ch-3$}
\put(2.0,.80){symplectic with $(\ch-c-2)$ basic classes}
\put(3.3,.63){(cf. \cite{FScan})}
\multiput(.170,.475)(.2,0){19}{$\bullet$}
\put (1.8, 0){Figure 15}
\end{picture}
\newpage

\noindent Simply connected complex surfaces are of three types:

\begin{itemize}
  \item homeomorphic to blowups of the complex projective plane $\CP$ or the Hirzebruch surfaces $F_n$, i.e. are rational surfaces,
  \item the elliptic surfaces and their logarithmic transforms $E(n)_{p,q}$, i.e. $c=0$, 
  \item are of general type, i.e. $2\ch-6 \le c < 9\ch$.
\end{itemize} We discussed the construction of the elliptic surfaces and their logarithmic transforms in Lectures 1 and 2 and discussed manifolds homeomorphic to rational surfaces in Lecture 5. We also showed how the operation of rational blowdown can give examples of irreducible simply connected symplectic manifolds with $0 < c <2\ch-6$ in Lectures 1 and 4. In fact every lattice point in the  $0 < c <2\ch-6$ region can be realized by an irreducible simply connected symplectic manifold \cites{FPS,FScan, stipeuro,parkgeo,stipbmy,GS}. The techniques used  are an artful application of either the generalized fiber sum or the rational blowdown constructions. Within this region there appears to be an interesting relationship between the minimal number of Seiberg-Witten basic classes and the pair $(\ch,c)$. In particular, all known smooth $4$-manifolds with $0  \le c \le \ch-3$ have at least $\ch - c -2 $ Seiberg-Witten basic classes \cite{FPS}. There is a physics, but yet no mathematical,  proof that this must be the case \cite{marino}. 

It is generally believed that all invariants in the range $2\ch-6 \le c \le 9\ch$ can be realized by a complex surface. It is a more delicate question to require these manifolds to be simply connected. It is known that all complex surfaces with $c= 9\ch$ and $\ch >1$ are ball quotients and hence cannot be simply connected. So the question remains if every lattice point in the $2\ch-6 \le c < 9\ch$ can be realized by a minimal (irreducible) simply connected complex manifold. There remain gaps, but this appears to be limited by the construction techniques rather than any intrinsic reason. In \cite{Parkgeo} Jongil Park uses the fiber sum and rational blowdown constructions to show that all  but at most finitely many lattice points in the $c \le 8\ch$ region can be realized by simply connected irreducible symplectic $4$-manifolds. For complex manifolds, the result is less sharp. For example Persson \cite{Perssongeo} shows that the region $2\ch-6 \le c \le 8(\ch - C\ch^{2/3})$ can be realized by a minimal simply connected complex surface where $C$ can be taken to be $9/12^{1/3}$.

The so called \emph{arctic region} where $8\ch < c < 9\ch$ has many gaps, again most likely due to a lack of sufficiently creative constructions. For simply connected complex surfaces see  \cite{persson}. For symplectic manifolds Andras Stipsicz \cite{stipbmy} and Jongil Park \cite{parkgeo} have constructions of simply connected symplectic $4$-manifolds in this region that start with nonsimply connected complex surfaces on or near the $c=9\ch$ line  and then use the fiber sum construction to kill the fundamental group. 

Later in this lecture we will concentrate on some interesting problems and constructions at the extremes of the region for surfaces of general type, i.e the Horikawa surfaces with $c=2\ch -6$ and manifolds with $c=9\ch$. Before we do there are some open problems that should not go unmentioned.

An irreducible $4$-manifold need not lie on a lattice point, for example just change the orientation of any manifold with odd $c$. Note that $\ch(X) \in \Z$ iff $X$ admits an almost-complex structure. In addition, the Seiberg-Witten invariants are only defined for manifolds with $\ch \in \Z$. Since our only technique to determine if a $4$-manifold is irreducible is to use  Seiberg-Witten invariants, all known irreducible $4$-manifolds have (with one orientation) $\ch \in \Z$.

Now irreducibility does not depend upon orientation. So the obvious question is whether there exist irreducible smooth 4-manifolds with $\ch \notin \Z$ for both orientations? Here the work of Bauer and Furuta \cite{BF} on stable homotopy invariants derived from the Seiberg-Witten equations may be useful. As a simple candidate for an irrecucible manifold with $\ch \notin \Z$ with both orientations consider two copies of the elliptic surface $E(2)$. Remove the neighborhood of a sphere with self-intersection $-2$ from each and glue together the resulting manifolds along their boundary ${{\mathbf{RP}}^{\,3}}$ using the orientation reversing diffeomorphism of ${{\mathbf{RP}}^{\,3}}$. The result has $\ch =\frac{7}{2}$ (or $\ch =\frac{37}{2}$ with the opposite orientation) and it is unknown if it is irreducible. See\cite{GS} for a handlebody description of this manifold.

As mentioned earlier, the work of Taubes \cite{TGW} shows that $c \ge 0$ for an irreducible symplectic $4$-manifold. Also, all of the known simply connected irreducible $4$-manifolds are homeomorphic (possibly changing orientation) to symplectic manifolds. There may be an interesting class of $4$-irreducible manifolds yet to be discovered, i.e. those with $c < 0$. (The above example has $c=-2$.) Note that if every simply connected irreducible smooth  $4$-manifolds has $c\ge 0$ regardless of orientation, then $c \le \frac{48}{5}\ch$. 

Figure 15 contains no information about the geography of spin $4$-manifolds, i.e. manifolds with $t=0$. For a spin $4$-manifold there is the relation $c = 8\ch \mod 16$. Almost every lattice point with $c = 8\ch \mod 16$ and $0 \le c < 9\ch$ can be be realized by an irreducible spin 4-manifold \cite{Park}. Surprisingly not all of the lattice points with $2\ch \le 3(\ch - 5)$ can be realized by complex manifolds with $t=0$ \cite{persson}, so spin manifolds with  $2\ch \le 3(\ch - 5)$ provide more examples of smooth irreducible $4$-manifolds with $2\ch-6 \le c < 9 \ch$ that support no complex structure (cf. \cites{Park, families}). 

\section{Uniqueness: The Botany Problem}

The uniqueness part of the classification scheme, sometimes referred to as {\it the botany problem},  is to determine all smooth (symplectic) $4$-manifolds with a fixed $(\ch,c, t)$ and to discover invariants that would distinguish them.  Here is where we begin to lose control of the classification of smooth $4$-manifolds.  As we have seen in our lecture on knot surgery, if  a topological $4$-manifold admits a smooth (symplectic) structure with nontrivial Seiberg-Witten invariants and that also contains a smoothly (symplectically) embedded minimal genus torus with self-intersection zero and with simply connected complement, then it also admits infinitely many distinct smooth (symplectic) structures and also admits infinitely many distinct smooth structures with no compatible symplectic structure. There are no known examples of (simply connected) smooth or symplectic $4$-manifolds with $\ch >1$ that do not admit such tori. Hence, there are no known smooth or symplectic $4$-manifolds with $\ch > 1$ that admit finitely many smooth or symplectic structures. We also showed in our last lecture that even for manifolds with $\ch =1$, for which there are no essential minimal genus tori with trivial normal bundle,  there are no known manifolds with finitely many smooth structures.

While knot surgery appears to be a new operation, we saw in an Lecture 3 that the knot surgery construction is actually a series of $\pm 1$ surgeries on {\it
nullhomologous tori}. The underlying observation is that any knot can be
unknotted via a sequence of crossing changes, which in turn can be
realized as a sequence of $\pm 1$ surgeries on unknotted curves
$\{c_1, \dots, c_n\}$ that link the knot algebraically zero times and
geometrically twice. When crossed with $S^{1}$ this translates to the
fact that $X$ can be obtained from $X_{K}$ via a sequence of $\pm 1$
surgeries on the {\it nullhomologous
tori} $\{S^1 \times c_1, \dots, S^1\times c_n\}$ in $X$. So the
hidden mechanism behind the knot surgery construction is surgeries on nullhomologous tori.  The calculation
of the Seiberg-Witten invariants is then reduced to understanding how
the Seiberg-Witten invariants change under a surgery on a nullhomologous torus. 

So if we expand the notion of torus surgery to include both homologically essential and nullhomologous tori, then we can eliminate the knot surgery construction from our list of the essential operations we have used to both construct smooth $4$-manifolds and to  alter their smooth structures. Thus our list is of essential operations is reduced to

\begin{itemize}
\item generalized fiber sum
\item surgery on a torus with trivial normal bundle
\item blowup
\item rational blowdown
\end{itemize}

There are further relationships between these operations. In
\cite{rat} it is shown that if $T$ is contained in a node
neighborhood, then a torus surgery can be
obtained via a sequence of blowups and rational blowdowns.  However, it is
not clear, and probably not true,  that a rational blowdown is always the result of blowups and
logarithmic transforms.

As an aside, this fact
together with work of Margaret Symington \cite{sym} shows that
logarithmic transformations ($p \ne 0$) on a symplectic torus results
in a symplectic manifold. We do not know of any general proof that a
logarithmic transformation on a symplectic torus in a
node neighborhood results in a symplectic manifold. However, there is an interesting proof that $E(1)_{p,q}$ is symplectic by observing that $E(1)_{p,q}$ is diffeomorphic to $E(1)_{T(p,q)}$ where $T(p,q)$ is the $(p,q)$ torus knot. The proof uses the fact that every diffeomorphism of the boundary of $E(1)\setminus F$ extends. A simple comparison of Seiberg-Witten invariants show that $E(n)_{p,q}$ is not diffeomorphic to $E(n)_{T(p,q)}$ for $n > 1$.

Rational blowdown changes the topology of the manifold $X$; while
$\ch$ remains the same, $c$ is decreased by $p-3$. So, an obvious
conjecture would be

\begin{conj}\label{log}
Two homeomorphic simply connected smooth 4-manifolds are related
via a sequence of surgeries on tori of square $0$.
\end{conj}

A slight weakening would be

\begin{conj}\label{logs}
Two homeomorphic simply connected smooth 4-manifolds are related
via a sequence of of surgeries on tori of square $0$ and rational blowups and rational blowdowns.
\end{conj}

Of course there are no counterexamples to these conjectures since these operations are the only ones used so far to change smooth structures. However, there are many examples of complex structures on a given manifold that are deformation inequivalent that are unknown whether they are diffeomorphic. A good test of these conjectures would be to determine if these complex structures are related via a sequence of surgeries on tori. 

\section{Horikawa Surfaces: How to go from one deformation type to another}

We mentioned in our first lecture that it is still unknown whether the two deformation types of Horikawa surfaces $H(4n-1)$ and $H^{'}(4n-1)$ with $(c,\ch,t)=(8n-8,4n-1)$ are, for $n>1$, diffeomorphic or even symplectomorphic. However, since both these deformation types are branched covers of $\SS$ we can use techniques of \cites{FSsurfaces, AKD} to braid the branched set to show that indeed $H(4n-1)$ can be obtained from $H^{'}(4n-1)$  via a $\pm 1$ surgery on an essential torus as follows.  

The basic idea is that the only difference between $H(4n-1)$ and $H^{'}(4n-1)$ is that both are covers of $\SS$ branched over a curve of bidgree $(6,4n)$, with the former a connected curve and the latter a disconnected curve. By braiding the later we make the curve connected. This braiding process is a procedure in $\SS$ which removes an (Lagrangian) annulus  whose boundary components lie on each component of the branch set for  $H^{'}(4n-1)$ and sews it back in so that the branch set is now connected. This annulus lifts to a torus in $H^{'}(4n-1)$ and the surgery operation lifts to a $\pm 1$ surgery on this torus. Of course, in order to know that the resulting manifold is $H^{'}(4n-1)$ we need to know that the resulting connected branch curve is isotopic to the connected complex curve of bidegree $(6,4n)$. This is precisely the situation covered by the work of Siebert and Tian \cite{SiebertTian} which guarantees that any two such symplectic curves in $\SS$ are isotopic.

 For a more general discussion of this phenomena and for other intriguing relationships between these deformation types see Auroux's fascinating article \cite{Aur}.

\section{Manifolds with $c=9\ch$: A fake projective plane}

One would like to state that every lattice point $(\ch, c)$ in the $c\le 9\ch$ can be realized by an irreducible simply connected $4$-manifold. This is almost true; however, there are open problems near the $c=9\ch$ line, called the {\it Bogomolov-Miyaoka-Yau (BMY)  line}. Part of the problem is that it is known that all complex surfaces with $c=9\ch > 9$ are ball quotients and hence not simply connected. As one approaches the $c=9\ch$ line it is not known that there exist complex surfaces, simply connected or not. Most of the constructions have the property that the closer you get to the $c=9\ch$ line the larger the fundamental group.  For more detailed information of this phenomena see \cite{persson} or the papers \cites{stipgeo, parkgeo, stipbmy, GS} that construct simply connected symplectic manifolds with large $c$ that asymptotically approach the $9\ch$ line. It would be very interesting to find simply connected symplectic manifolds with $c=9\ch > 9$ or with $c > 9\ch$.

Part of the issue here is that the all the constructions of manifolds with positive signature, except for $\CP$,  rely on branched covering constructions with singular branched set and the closer one gets to the $c=9\ch$ line the more singular the branched set. These singularities increase the complexity of the fundamental group of the cover. In the remainder of this lecture we would like to give the simplest construction that we know of a fake manifold with $c=9\ch$. This construction is due to Keum \cite{keum} and starts with the existence of a particular configuration of spheres on a Doglachev surface shown to exist by Ishida \cite{ishida}.  The result will be a complex surface with the integral homology of $\CP$. 

In \cite{ishida} Ishida constructs a Dolgachev surface $E(1;2,3)$ with a remarkable property. This elliptic fibration has two multiple fibers of multiplicity $2$ and $3$ respectively, and has an elliptic fibration with four singularities of type $I_3$ (see Lecture 5). It is not difficult to show that $E(1;2,3)$ has a sextuple section which is a sphere of self-intersection $-3$. The remarkable property of Ishida's manifold is that it contains the following four disjoint configuration of spheres: 

\centerline{\unitlength 1cm
\begin{picture}(5,4.5)
\put(1.1,.7){$\bullet$}
\put(1.2,.8){\line(1,0){1.64}}
\put(2,.7){$\bullet$}
\put(2.8,.7){$\bullet$}
\put(1,1.1){$-3$}
\put(1.75,1.1){$-2$}
\put(2.55,1.1){$-2$}
\put(1.4,1.7){$\bullet$}
\put(1.5,1.8){\line(1,0){.75}}
\put(2.25,1.7){$\bullet$}
\put(1.3,2){$-2$}
\put(2.1,2){$-2$}
\put(1.4,2.7){$\bullet$}
\put(1.5,2.8){\line(1,0){.75}}
\put(2.25,2.7){$\bullet$}
\put(1.3,3){$-2$}
\put(2.1,3){$-2$}
\put(1.4,3.7){$\bullet$}
\put(1.5,3.8){\line(1,0){.75}}
\put(2.25,3.7){$\bullet$}
\put(1.3,4){$-2$}
\put(2.1,4){$-2$}
\end{picture}}

\noindent The first three configurations sit inside three of the four $I_3$ singularities and the $-3$ sphere in the last configuration is the sextuple section with the $-2$ spheres coming from the fourth $I_3$ singularity.  It is remarkable to find such a configuration of nine spheres inside the Dolgachev surface $E(1;2,3)$ since they capture all its $b^{-}$. 

With this starting point Keum then observes that the complement of the first three configurations has euler characteristic $3$ and that it has a $3$-fold cover with three boundary components that are $S^3$ each of which can be capped with a $4$-ball to obtain a manifold $X$ with euler characteristic $12$. He shows that this is again a Doglachev surface with four singular fibers, three of type $I_1$ and one of type $I_9$. The important feature is that $X$ now contains the following three disjoint configuration of spheres:

\centerline{\unitlength 1cm
\begin{picture}(5,4)
\put(1.1,.7){$\bullet$}
\put(1.2,.8){\line(1,0){1.7}}
\put(2,.7){$\bullet$}
\put(2.8,.7){$\bullet$}
\put(1,1.1){$-3$}
\put(1.75,1.1){$-2$}
\put(2.55,1.1){$-2$}
\put(1.1,1.7){$\bullet$}
\put(1.2,1.8){\line(1,0){1.7}}
\put(2,1.7){$\bullet$}
\put(2.8,1.7){$\bullet$}
\put(1,2.1){$-3$}
\put(1.75,2.1){$-2$}
\put(2.55,2.1){$-2$}
\put(1.1,2.7){$\bullet$}
\put(1.2,2.8){\line(1,0){1.7}}
\put(2,2.7){$\bullet$}
\put(2.8,2.7){$\bullet$}
\put(1,3.1){$-3$}
\put(1.75,3.1){$-2$}
\put(2.55,3.1){$-2$}
\end{picture}}

\noindent each a lift from $E(1;2,3)$. Again, this is a remarkable configuration to find in a Doglachev surface since it again captures all of $b^{-}$. To finish Keum's construction he observes that the complement of this configuration has euler characteristic zero and has a $7$-fold cover  with three boundary components each a $3$-sphere $S^3$ which can again be capped off with a $4$-ball to obtain a manifold $Y$ with Euler characteristic $3$ and with $b_1=0$. This is Keum's construction of a fake (homology) projective plane. Of course there is some heavy lifting in showing the existence  of the Ishida manifold for which we do not know of a direct (topological) construction.

\section{Small $4$-manifolds} As we have already mentioned, finding simply connected manifolds close to the $c=9\ch$ line remains a challenge. This challenge is exposed in our inability to find irreducible simply connected $4$-manifolds with small Euler characteristic. We have discussed in Lecture 5 the case when $\ch =1$, i.e. are there irreducible simply connected smooth $4$-manifolds homeomorphic to $n$ blowups of $\CP$ with $n<5$? For 
$\ch=2$ the best that one can do to date is $c=11$ \cites{Park3, SSb}, i.e. there are irreducible simply connected smooth $4$-manifolds homeomorphic to $3\CP \# n\CPb$ for $8 \le n \le 19$. The constructions are similar to those discussed in Lecture 5 applied to $E(2)$.  What about $n <8$?. The inability to get close to the $c=9\ch$ line permeates for  $\ch >2$. Further progress on constructing simply connected $4$-manifolds with small Euler characteristic should yield further insight into constructing manifolds in the arctic region.

\section{What were the four $4$-manifolds?}

Our four $4$-manifolds represent examples of simply connected complex surfaces with Kodaira dimension $-\infty$, $0$, $1$, and $2$, {\em{i.e.}} $\CP$ and its blowups, the K3 surface $E(2)$, the elliptic fibrations $E(n)$, $n >2$, and $H(4n-1)$, $n>1$.

\end{document}